\documentclass[10pt,a4paper]{article}
\usepackage{amstext}
\usepackage{array}
\usepackage{amsfonts}
\usepackage{amssymb, bm}
\usepackage{graphicx}
\usepackage{epstopdf}
\usepackage{algorithm}     
\usepackage{algpseudocode} 
\usepackage{url}
\usepackage{caption}
\usepackage{subfig}
\usepackage{csquotes}
\usepackage{tikz}
\usetikzlibrary{arrows}
\usetikzlibrary{decorations.markings}
\usetikzlibrary{shapes.geometric}
\usetikzlibrary{shapes.arrows}
\usetikzlibrary{fit,calc}
\usetikzlibrary{quotes,angles}
\usetikzlibrary{positioning, arrows.meta}
 \usetikzlibrary{patterns}
\usepackage{animate}  
\usepackage{graphicx}
\usepackage{pgf}
\usepackage{bbm}
\usepackage{pgfplots}
\usepackage{mathtools}
\usepackage{physics}
\usepackage{amsthm}
\usepackage{hyperref}
\usepackage{stmaryrd}

\DeclareMathOperator{\card}{\#}
\newcommand{\oo}{\overline\Omega}

\newcommand{\vertiii}[1]{{\left\vert\kern-0.25ex\left\vert\kern-0.25ex\left\vert #1 
   \right\vert\kern-0.25ex\right\vert\kern-0.25ex\right\vert}}

\newtheorem{Definition}{Definition}[section]
\newtheorem{theorem}{Theorem}[section]
\newtheorem{theorem1}{Theorem}[section]
\newtheorem{theorem2}{Theorem}[section]
\newtheorem{proposition}[theorem]{Proposition}
\newtheorem{lemma}[theorem1]{Lemma}
\newtheorem{remark}[theorem2]{Remark}

\usepackage{geometry}
 \usepackage{fancyhdr}  
  \fancypagestyle{plain}{
\fancyhead{}
\fancyhead[C]{\hfill Submitted to SIAM Journal on Numerical Analysis, January 2026}
 }

\title{Mathematical aspects of registration methods in bounded domains}
\author{Angelo Iollo,  Jon  Labatut, Pierre Mounoud, Tommaso Taddei}
\date{\today}

\begin{document}

\maketitle

\begin{abstract}
Registration methods in bounded domains have received significant attention 
in the model reduction literature,
as a valuable tool for nonlinear approximation.
The aim of this work is to provide a concise yet complete overview of relevant results  for registration methods in  $n$-dimensional domains,  from the perspective of parametric model reduction.  
We present a thorough analysis of two classes of methods, vector flows and compositional maps: we discuss the enforcement of the bijectivity constraint and we comment on the approximation properties of the two methods, for Lipschitz  $n$-dimensional domains.
\end{abstract}

\section{Introduction}
\label{sec:intro}
In the past few decades, 
registration methods have received considerable attention in scientific computing. 
In computer vision and pattern recognition, registration refers to the process of finding a spatial transformation that aligns two datasets \cite{brown1992survey,ma2015robust,myronenko2010point}.
Registration techniques share important features 
with   morphing techniques in mesh adaptation
\cite{mcrae2018optimal} and shock-tracking methods in CFD
\cite{shubin1982steady,zahr2018optimization}.
More recently, 
registration methods have gained significant attention in the field of  (parametric) model order reduction (MOR): in this context, registration techniques have been proposed to find a parametric transformation that improves
the compressibility of the solution set associated with a  given parametric partial differential equation (PDE)
\cite{blickhan2023registration,cucchiara2024model,iollo2022mapping,mojgani2021low,taddei2020registration}.
Registration techniques for MOR applications should
accurately preserve   the domain geometry to  enable
the proper enforcement of boundary conditions; furthermore, they should meet the computational and memory requirements of the \emph{offline/online} computational decomposition.

We denote by  $\Omega$ a Lipschitz bounded domain in $\mathbb{R}^n$.
We introduce the  set of diffeomorphisms  
${\rm Diff}(\overline{\Omega})$:
 the function $\Phi: \overline{\Omega} \to \mathbb{R}^n$ belongs to ${\rm Diff}(\overline{\Omega})$ if it is of class $C^1$ up to the boundary and is bijective  
  (one-to-one and onto) in $\overline{\Omega}$.
  For MOR applications, given the parameter space $\mathcal{P} \subset \mathbb{R}^{n_\mu}$, we are interested in parametric diffeomorphisms 
  $\Phi: \overline{\Omega} \times \mathcal{P} \to \overline{\Omega}$, which satisfy 
  (i) $\Phi(\cdot, \mu)\in {\rm Diff}(\overline{\Omega})$ for all $\mu\in \mathcal{P}$,  
  and 
  (ii) $\Phi$ is continuous with respect to the parameter.
If we introduce the parametric target (or matching) function 
 $\mathfrak{f}^{\rm tg}: C(\overline{\Omega}; \mathbb{R}^n) \times \mathcal{P} \to \mathbb{R}$, the goal of \emph{optimization-based} \cite{taddei2025compositional} or
 \emph{variational} \cite{beg2005computing} methods for registration is to find a parametric bijection that approximately minimizes
\begin{equation}
\label{eq:optimization_based_registration_ideal}
\min_{\Phi \in {\rm Diff}(\overline{\Omega})} 
\mathfrak{f}^{\rm tg}(\Phi; \mu),
\end{equation}
for all parameters $\mu\in \mathcal{P}$.
For MOR, $\mathfrak{f}^{\rm tg}$ should be informed by the  relevant features of the solution field we wish to track; in image registration, $\mathfrak{f}^{\rm tg}(\Phi; \mu)$ measures the discrepancy between the target object and the deformed reference object --- we provide two practical examples of 
target functions in section \ref{sec:ansatz_registration}.

Since the set $ {\rm Diff}(\overline{\Omega})$ is infinite-dimensional, the point of departure of registration techniques is to identify an effective finite-dimensional parameterization that is amenable for computions. 
Towards this end, we introduce
 the function $\texttt{N}: \mathbb{R}^M \to {\rm Lip}(\overline{\Omega}; \mathbb{R}^2)$ and a  penalty function $\mathfrak{f}_{\rm pen}:\mathbb{R}^M \to \mathbb{R}_+$ such that
\begin{equation}
\label{eq:desiderata_mapping_space}
\left\{
\begin{array}{l}
\displaystyle{
\mathcal{B}_{\texttt{N}}
=\left\{
\texttt{N}(\mathbf{a})  \; :  \;
\mathbf{a} \in \mathbb{R}^M, \; 
\mathfrak{f}_{\rm pen}(\mathbf{a}) \leq C
\right\}
\subset {\rm Bij}(\overline{\Omega}),
\quad
{\rm for \; some} \; C>0;
}
\\[3mm]
\displaystyle{
\texttt{N}(\mathbf{a}= 0) = \texttt{id},
\;\;
\mathfrak{f}_{\rm pen}(\mathbf{a}= 0)  <  C;
}
\\[3mm]
\displaystyle{
\texttt{N}, \; \mathfrak{f}_{\rm pen}
\; {\rm are \; Lipschitz \; continuous}.
}
\\
\end{array}
\right.
\end{equation}
Here, ${\rm Bij}(\overline{\Omega})$ refers to the set of bijections in $\overline{\Omega}$.
Then, we introduce the surrogate of 
 \eqref{eq:optimization_based_registration_ideal} 
\begin{equation}
\label{eq:tractable_optimization_based_registration_morozov}
\min_{   \mathbf{a}  \in \mathbb{R}^M}
 \mathfrak{f}_{\mu}^{\rm tg}(\texttt{N}( \mathbf{a} ))
\; {\rm s.t.} \;\; 
 \mathfrak{f}_{\rm pen}(\mathbf{a} ) \leq C,
\end{equation}
which can be solved using standard optimization techniques (e.g., interior point methods).
To further simplify the optimization task, we might also consider the Tikhonov regularization of \eqref{eq:tractable_optimization_based_registration_morozov}, 
\begin{equation}
\label{eq:tractable_optimization_based_registration}
\min_{   \mathbf{a}  \in \mathbb{R}^M}
 \mathfrak{f}_{\mu}^{\rm obj}(\mathbf{a} ):=
 \mathfrak{f}_{\mu}^{\rm tg}(\texttt{N}( \mathbf{a} ))
\; +\lambda  \; 
 \mathfrak{f}_{\rm pen}(\mathbf{a} ),
\end{equation}
where the hyper-parameter $\lambda >0$  should be properly tuned to ensure the bijectivity constraint.

Different registration techniques differ due to the choice of the pair
$(\texttt{N}, \mathfrak{f}_{\rm pen})$.
A first set of registration strategies 
is based on the parameterization of a pseudo-velocity field $v=v(\mathbf{a})$ and on the definition of $\texttt{N}(\mathbf{a})$ as the flow of the vector field $v(\mathbf{a})$.
A second set of strategies exploits a direct parameterization of the displacement field: to cope with 
curved domains, these strategies introduce a transformation from $\Omega$ to a domain 
$\Omega_{\rm p}$ with piecewise-linear boundaries.
Below, we refer to the first class of maps as
 \emph{vector flows} (VFs), and to the second class as \emph{compositional maps} (CMs).
 
 The development of 
 VFs exploits the recent progresses in the field of optimal transport \cite{peyre2019computational,santambrogio2015optimal}.
Since the seminal works of Christensen 
 \cite{christensen1996deformable} and Beg
\cite{beg2005computing}, 
 VFs have been widely adopted  for image registration tasks \cite{sotiras2013deformable}; 
 VF maps have also been considered for
mesh adaptation  \cite{de2016optimization}, geometry  reduction 
\cite{kabalan2025elasticity}, and
  model reduction in \cite{iollo2014advection,labatut2025non}. In particular, the work of Labatut \cite{labatut2025non} explicitly addresses the issue of adequately representing the boundary of the domain, for two- and three-dimensional geometries.
On the other hand, CMs have been introduced in the MOR literature in
  \cite{taddei2020registration,taddei2025compositional}, and further developed  and adopted in \cite{mirhoseini2023model,razavi2025registration}.

The objective of this paper is to present  key mathematical properties of VFs and CMs for registration,  from the perspective of parametric model reduction.  
The ultimate goal is  a mathematically-sound comparison of state-of-the-art techniques for parametric registration.
We also aim to identify key technical limitations of available techniques that might inspire future contributions to  the subject.
Finally, we investigate the properties of several modal expansions for both VFs and CMs: the goal here is to reduce \emph{a priori} the dimensionality of the registration problem.

The analysis of VFs  exploits well-known results (see, e.g.,  \cite{younes2010shapes}): 
in contrast to  \cite{younes2010shapes},  we here state the results in the finite-dimensional setting and we address in detail the issue of approximation (cf. Proposition \ref{th:approx_flows}).
On the other hand, the results for CMs are new and  
extend the analysis of \cite{taddei2025compositional} to $n$-dimensional geometries.
{These results provide the rigorous foundation for a thorough comparison of the two classes of methods.}
Finally, the discussion on the choice of the basis is original.

The outline of the work is as follows.
Section \ref{sec:ansatz_registration} introduces VFs and CMs for registration, and provides two practical examples of target functions;
section \ref{sec:analysis} states the theoretical results of this work and provides a detailed discussion of the two strategies.
Section \ref{sec:truncated_basis} discusses two methods for the generation of the modal expansion for both VFs and CMs, and
 provides  numerical experiments to 
 assess
  the accuracy of the modal bases
 for a representative test case.
Section \ref{sec:proofs} provides the proofs of the results of section \ref{sec:analysis}. The  appendix provides additional proofs that either 
appeared in different contexts in previous works
or are immediate consequences of well-known results.

\section{Formulation}
\label{sec:ansatz_registration}
In this section, we discuss two strategies to effectively parameterize diffeomorphisms in bounded domains: \emph{vector flows (VFs)} and \emph{compositional  maps (CMs)}. Given the Lipschitz connected domain $\Omega\subset \mathbb{R}^n$, we denote by $\mathbf{n}: \partial \Omega \to \mathbb{S}^2 = \{x \in \mathbb{R}^n: \|  x \|_2=1  \}$ the outward unit normal, and by $C^k(\overline{\Omega}; \mathbb{R}^n)$ the space of $k$-times differentiable functions up to the boundary, with values in $\mathbb{R}^n$. 
We denote by $\mathcal{U}_0$ the subspace of $C^1$ vector fields  with normal displacement equal to zero,
 $$
 \mathcal{U}_0:= \left\{
 \varphi \in C^1(\overline{\Omega}; \mathbb{R}^n) \, : \,
 \varphi \cdot \mathbf{n} \big|_{\partial \Omega}  = 0
 \right\},
 $$
 we further introduce the space of time-dependent fields
 $$
 \mathcal{V}_0:=  \{ v\in C^1(\overline{\Omega} \times [0,1]; \mathbb{R}^n) \,: \, v(\cdot,t) \cdot \mathbf{n} \big|_{\partial \Omega} = 0 \;\;  \forall \, t\in [0,1] \}
 $$ 
Below, $\texttt{id}:\Omega \to \Omega$ denotes the identity map,
$\texttt{id}(x)=x$, while $\mathbbm{1}\in \mathbb{R}^{n\times n}$ is the identity matrix.
For Lipschitz domains, the normal vector is 
defined almost everywhere (a.e.)
(cf. Rademacher's theorem  \cite{ziemer2012weakly}); therefore, the conditions 
 $\varphi \cdot \mathbf{n} \big|_{\partial \Omega} =0 $ and 
$v(\cdot, t) \cdot \mathbf{n} \big|_{\partial \Omega} =0 $ should also be interpreted in an a.e. sense.
Finally, we introduce an important class of diffeomorphisms in $\mathbb{R}^n$ that are important for the discussion.

\begin{Definition}
\label{def:Diff0}
(\cite{banyaga2013structure})
Let $\Phi$ be a diffeomorphism in the Lipschitz domain $\Omega\subset \mathbb{R}^n$. We say that $\Phi$ is isotopic to the identity if there exists a smooth field (isotopy) $X: \overline{\Omega} \times[0,1] \to  \mathbb{R}^n$ such that
(i) \emph{$X(\cdot, 0)=\texttt{id}$}, 
(ii) $X(\cdot, 1)=\Phi$, 
(iii) the map $X(\cdot, t) \in {\rm Diff}(\overline{\Omega})$ for all $t\in [0,1]$;
(iv) the map $\partial_t X(\cdot, t) \in C^1(\overline{\Omega})$ for all $t\in [0,1]$.
We denote by ${\rm Diff}_0(\overline{\Omega})$ the subset of diffeomorphisms that are isotopic to the identity.  
\end{Definition}

As discussed below (cf. Propositions \ref{th:evolution_jacobian} and \ref{th:approximation_DB}), we can easily construct diffeomorphisms that do not belong to ${\rm Diff}_0(\overline{\Omega})$;  however, 
${\rm Diff}_0(\overline{\Omega})$ constitutes the relevant set of diffeomorphisms for model reduction applications.
Given the   parameter domain $\mathcal{P}\subset \mathbb{R}^{n_\mu}$, 
in model order reduction, we seek parametric mappings
$\Phi: \overline{\Omega} \times \mathcal{P} \to \overline{\Omega}$ such that
$\Phi(\cdot,\mu) :\overline{\Omega}  \to \overline{\Omega}$ is a diffeomorphism for all $\mu\in \mathcal{P}$
\cite{taddei2020registration}.
Exploiting Definition \ref{def:Diff0}, we can readily show
that
 if 
(i) 
 $\mathcal{P}$ is connected,
 (ii) $\Phi(\cdot,\mu')$ is equal to the identity for some $\mu'\in \mathcal{P}$, and
 (iii) $\Phi$ is continuous with respect to the parameter,  we must have that 
 $\Phi(\cdot,\mu)\in {\rm Diff}_0(\overline{\Omega})$  for all $\mu\in \mathcal{P}$.
  This observation shows the practical relevance of  this subset of diffeomorphisms for model reduction applications.

In this work, we consider diffeomorphisms up to the boundary of $\Omega$: in scientific computing, it is indeed important to ensure that the boundary of $\Omega$ is mapped in itself to properly enforce boundary conditions. Note that we can readily construct diffeomorphisms in  ${\rm Diff}_0(\Omega)$ that do not admit an extension to $\overline{\Omega}$  as shown in the next remark.

\begin{remark}
\label{remark:counterexample}
We provide an example of diffeomorphism in  ${\rm Diff}_0(\Omega)$ that does not admit an extension to $\overline{\Omega}$.  
Consider $\Omega=\{x\in \mathbb{R}^2:  \|x\|_2 < 1\}$ and define the  function 
$\Phi:\Omega \to \Omega$ such that
$$
\Phi(\xi)  = r \left[
\begin{array}{l}
    \cos( \theta + \phi(r) ) \\
   \sin( \theta + \phi(r) ) \\
\end{array}
\right]
\quad
{\rm with} \;
r=\|\xi\|_2, \;\;
\theta = \arctan_2(\xi),
$$
and $\arctan_2: \mathbb{R}^2 \to (-\pi,\pi)$ is the two-argument arctangent, and $\phi:[0,1)$ is a $C^\infty(0,1)$ function such that $\phi(r)=0$ for $r<1/2$ and $\lim_{r\to 1}\phi(r) = +\infty$. Since $\phi\equiv 0$ in a neighborhood of $r=0$, we find that $\Phi$ is the identity in a 
neighborhood of $\xi=0$ and thus it
is of class $C^1$ in $\Omega$. We also observe that  $\Phi$ is  a bijection from  each circle of radius $r$ in itself; therefore,  $\Phi \in {\rm Diff}(\Omega)$. In addition, 
$$
X(\xi,t) = r \left[
\begin{array}{l}
    \cos( \theta + t \phi(r) ) \\
   \sin( \theta + t \phi(r) ) \\
\end{array}
\right],
$$
is an isotopy from the identity to $\Phi$; we conclude that $\Phi \in {\rm Diff}_0(\Omega)$.
However,  $\Phi$ does not admit a limit for any $x\in \partial \Omega$: therefore, it does not admit  a continuous extension to $\overline{\Omega}$.   
\end{remark}

\subsection{Vector flows for registration}
\label{sec:velocity_based_map}
Let 
 $v:\overline{\Omega} \times [0,1]\to \mathbb{R}^n$ 
satisfy  $v \in C(\overline{\Omega} \times [0,1]; \mathbb{R}^n)$, $v(\cdot,t) \in C^1(\overline{\Omega}; \mathbb{R}^n)$ for all $t\in [0,1]$, and
$v(\cdot,t) \cdot \mathbf{n} \big|_{\partial \Omega} = 0$.
Then, we say that $X:\overline{\Omega}\times [0,1]\to \overline{\Omega}$ is the flow of the vector field $v$ if 
\begin{equation}
\label{eq:flow_diffeomorphisms}
\left\{
\begin{array}{ll}
\frac{\partial X}{\partial t}(\xi,t) = v( X(\xi,t) , t) & t\in (0,1], \\[3mm]
X(\xi,0) = \xi, & \\
\end{array}
\right. 
\quad
\forall \, \xi \in \overline{\Omega}.
 \end{equation}
 Furthermore, we define 
$F[v]: \overline{\Omega} \to \overline{\Omega}$ such that 
 $F[v](\xi):=X(\xi,t=1)$. In this approach, we seek diffeomorphisms of the form:
 \begin{equation}
 \label{eq:VB_maps}
 \texttt{N}(\xi; \mathbf{a}) := F[v(\mathbf{a})](\xi), \quad
 {\rm where} \;\;
 v(x, t; \mathbf{a})
 =\sum_{i=1}^M \, (\mathbf{a})_i \phi_i(x,t),
 \end{equation}
 where $\mathbf{a}\in \mathbb{R}^M$ is a vector of coefficients,
 and 
 $\phi_1 \ldots,\phi_M \in  \mathcal{V}_0(\Omega)$ 
 are linearly-independent.
 
\subsection{Compositional maps for registration}
\label{sec:displacement_based_map}
As opposed to VFs, CMs rely on a parameterization of the displacement field.
We distinguish between two scenarios.
If $\Omega$ is a polytope (i.e., a polygon if $n=2$ or a polyhedron if $n=3$), we define
\begin{equation}
\label{eq:DB_ansatz_polytope}
\texttt{N}( \xi;   \mathbf{a}) = \xi + \sum_{i=1}^M (\mathbf{a})_i \varphi_i(\xi),
\end{equation}
where
  $\varphi_1 \ldots,\varphi_M \in \mathcal{U}_0(\Omega)$ are linearly-independent. For general curved domains, 
  we introduce the polytope
$\Omega_{\rm p}$ isomorphic to $\Omega$ and the Lipschitz bijection $\Psi: \Omega_{\rm p} \to \Omega$. Then, we define the CM ansatz:
\begin{equation}
\label{eq:DB_ansatz}
\texttt{N}(\mathbf{a}) = \Psi \circ  \texttt{N}_{\rm p}(\mathbf{a}) \circ \Psi^{-1},
\qquad
\texttt{N}_{\rm p}(\xi; \mathbf{a}) = \xi + \sum_{i=1}^M (\mathbf{a})_i \varphi_i(\xi),
\end{equation}
where
  $\varphi_1 \ldots,\varphi_M \in \mathcal{U}_0(\Omega_{\rm p})$ are linearly-independent.
We refer to \cite{taddei2025compositional} for a thorough discussion on the construction of the mapping $\Psi$ and the polytope $\Omega_{\rm  p}$ for two-dimensional domains: the extension to three-dimensional domains is the subject of ongoing research.

\subsection{Representative target functions}
We provide two practical examples of target functions used in both image registration (see, e.g., \cite{younes2010shapes}) and model reduction \cite{taddei2025compositional}. 
 In more detail, the first target function is used in
 \cite{taddei2020registration,taddei2025compositional} in the context of model reduction and
 generalizes the target function used in \cite{beg2005computing} for image registration, while the second target is broadly used for point set registration \cite{myronenko2010point}.
 
Given the field $u\in H^1(\Omega; \mathbb{R})$ and the $N$-dimensional linear space 
$\mathcal{Z}_N \subset  L^2(\Omega; \mathbb{R})$, we 
 introduce the \emph{distributed objective}: 
 \begin{equation}
 \label{eq:distributed_sensor}
 \mathfrak{f}_{(1)}^{\rm tg}(\Phi) = 
 \min_{\zeta\in \mathcal{Z}_N}
\frac{1}{2} 
 \int_\Omega ( u\circ \Phi(\xi) -\zeta(\xi) )^2 \, d\xi.
 \end{equation}
 The linear space $\mathcal{Z}_N$ encodes our knowledge of the target and can be built using the adaptive algorithm introduced in \cite{taddei2021space}.

Given the set of points $\{  \xi_i \}_{i=1}^{N_0}$ and
$\{ y_j \}_{j=1}^{N_1}$, and the matrix of weights 
$\mathbf{P}  = [P_{i,j}]_{i,j} \in \mathbb{R}^{N_1\times N_0}$ such that
$P_{i,j}\in [0,1]$, and
$\sum_{j} P_{i,j}\ =1$ and 
$\sum_{i} P_{i,j} = 1$ for all $i=1,\ldots,N_0$ and $j=1,\ldots,N_1$, 
we introduce the 
\emph{pointwise objective}:
 \begin{equation}
 \label{eq:pointwise_sensor}
 \mathfrak{f}_{(2)}^{\rm tg}(\Phi) = 
 \frac{1}{2} 
\sum_{i=1}^{N_0}  \sum_{j=1}^{N_1}
P_{i,j} \, \big\|
 \Phi(\xi_i) - y_j
\big\|_2^2;
 \end{equation}
 as discussed in \cite{myronenko2010point} the weights
 $\{ P_{i,j} \}_{i,j}$ describe the likelihood  of the point $\xi_i$ to be mapped in the point $y_j$ and are adaptively selected using an expectation-maximization (EM) procedure.
 
\section{Mathematical analysis}
\label{sec:analysis}
\subsection{Existence of minimizers}

We introduce the    Banach space $V \subset C(\overline{\Omega}; \mathbb{R}^n)$ with norm $\|\cdot \|_V $.
We assume that 
 the restrictions   $\phi_1(\cdot, t=1),\ldots,\phi_M(\cdot, t=1) $ of the fields in \eqref{eq:VB_maps} belong to   $V$;
 similarly, we assume that 
 $\varphi_1,\ldots,\varphi_M$ in  \eqref{eq:DB_ansatz_polytope} belong to   $V$.  
 Finally, we assume that the target function $\mathfrak{f}^{\rm tg}$ 
in \eqref{eq:optimization_based_registration_ideal} 
 is Fréchet-differentiable with respect to $V$.
We recall 
 that $\mathfrak{f}^{\rm tg}$ is Fréchet-differentiable in $V$ if, 
 for all $\Phi\in V$, there exists a linear continuous  operator
  $D \mathfrak{f}^{\rm tg}[\Phi]: V \to \mathbb{R}$ such that
 \begin{equation}
 \label{eq:frechet_derivative}
\lim_{\| h \|_V \to 0} \frac{\big|   \mathfrak{f}^{\rm tg}[\Phi+ h] -        
 \mathfrak{f}^{\rm tg}[\Phi] -  
 D \mathfrak{f}^{\rm tg}[\Phi]  (h)    \big|}{\| h \|_V}  = 0.
 \end{equation}
 
Note that any Fr{\'e}chet differentiable function is continuous; furthermore, if $\mathbf{a} \in \mathbb{R}^M \mapsto \texttt{N}(\mathbf{a}) \in V$ is differentiable, we can apply the chain rule to deduce an expression for the derivative of $E: \mathbf{a} \mapsto \mathfrak{f}^{\rm tg} ( \texttt{N}(\mathbf{a}) )$, 
 \begin{equation}
\label{eq:sensitivity_method_b}
\frac{\partial E}{\partial a_i} (\mathbf{a})
\; = \; 
 D \mathfrak{f}^{\rm tg} 
 \left[  \texttt{N}(\mathbf{a}) \right] 
\left( 
  \frac{\partial \texttt{N}}{\partial a_i} 
(\mathbf{a})
\right),
\quad
i=1,\ldots,M.
\end{equation}

It is easy to verify that \eqref{eq:distributed_sensor} and \eqref{eq:pointwise_sensor} are Fr{\'e}chet differentiable in $V = {\rm Lip}(\Omega; \mathbb{R}^n)$. 
By tedious but straightforward calculations, we can indeed find the following explicit expressions of the Fréchet derivative
of \eqref{eq:distributed_sensor} and \eqref{eq:pointwise_sensor}:
\begin{equation}
 \label{eq:sensors_der}
\begin{array}{l}
\displaystyle{
D \mathfrak{f}_{(1)}^{\rm tg}[\Phi]  ( h)
 = 
 \int_\Omega ( u\circ \Phi(\xi) -
\zeta_\Phi 
  ) 
\nabla u \circ \Phi(\xi) h(\xi) \,  
  \, d\xi,
} \\[3mm]
\displaystyle{
  D \mathfrak{f}_{(2)}^{\rm tg}[\Phi]  \,    (h)
\,  = \,
\sum_{i=1}^{N_0} 
\left(
 \Phi(\xi_i)   \, \, -
 \sum_{j=1}^{N_1}
P_{i,j}
 \,  y_j
\right)
\cdot  h(\xi_i) ,
}
\\
\end{array} 
 \end{equation}
where $\zeta_\Phi : =  {\rm arg} \min_{\zeta\in \mathcal{Z}_N}
 \int_\Omega ( u\circ \Phi(\xi) -\zeta(\xi) ) 
  \, d\xi$ is the $L^2$ projection  of 
$ u\circ \Phi$ onto $\mathcal{Z}_N$.

Next result shows that the problem 
\eqref{eq:tractable_optimization_based_registration} admits a minimizer $\mathbf{a}^\star\in \mathbb{R}^M$ under reasonable assumption on the penalty  $\mathfrak{f}_{\rm pen}$ and the target function $\mathfrak{f}^{\rm tg}$. Note that the solution to
\eqref{eq:tractable_optimization_based_registration} 
is not guaranteed to be unique.

\begin{proposition}
\label{th:existence_minimizers}
Suppose that \emph{$\texttt{N}$} and 
$\mathfrak{f}_{\rm pen}$ satisfy \eqref{eq:desiderata_mapping_space}. Suppose also that
(i) 
\emph{$E:= \mathbf{a} \mapsto \mathfrak{f}^{\rm tg}(\texttt{N}(\mathbf{a}))$} is Lipschitz continuous and bounded from below, and
(ii) $\displaystyle{\lim_{\| \mathbf{a} \|_2 \to \infty}}$
$\mathfrak{f}_{\rm pen}(\mathbf{a}) = +\infty$. Then, the function $\mathbf{a} \mapsto  \mathfrak{f}^{\rm obj}(\mathbf{a})$  in 
\eqref{eq:tractable_optimization_based_registration}  has a minimum in $\mathbb{R}^M$.
\end{proposition}

\begin{proof}
The result is a direct consequence of Weierstrass theorem; see  section \ref{sec:easy_results}.
\end{proof}

\subsubsection{Application of gradient descent procedures}
Exploiting the discussion above, we might be tempted to directly use a gradient-descent technique   to    find  a local  minimum
of \eqref{eq:tractable_optimization_based_registration}
 (we omit the penalty term to shorten notation),
\begin{equation}
\label{eq:wrong_gradient_descent}
\mathbf{a}^{(k+1)}
=
\mathbf{a}^{(k)}
- \gamma^{(k)}
\,
\nabla E(\mathbf{a}^{(k)}),
\quad
k=1,2,\ldots
\end{equation}
for a proper choice of the step size $\gamma^{(k)} \in (0,1]$.
We empirically find that the update \eqref{eq:wrong_gradient_descent} performs poorly in practice:
we should hence find an appropriate preconditioner to speed up convergence.
Towards this end, we  exploit the variational formulation of \eqref{eq:tractable_optimization_based_registration}.
Below, we only consider VFs; however, the same reasoning applies to CMs.

We introduce the reduced space $\mathcal{W}_M  = {\rm span} \{ \phi_i \}_{i=1}^M$ and the Hilbert space $H\subset \mathcal{V}_0(\Omega)$ with
inner product $(\cdot, \cdot)_H$ and duality pairing  $\langle \cdot, \cdot \rangle_H$. Given the velocity fields $v,h\in \mathcal{W}_M$, we denote by $\mathbf{a}, \mathbf{h}\in \mathbb{R}^M$ the corresponding vectors of coefficients ---
that is,
$v = \sum_i (\mathbf{a})_i \phi_i$ and 
$h = \sum_i (\mathbf{h})_i \phi_i$; next, we restate \eqref{eq:tractable_optimization_based_registration} as
\begin{equation}
\label{eq:optimization_variational}
\min_{v\in \mathcal{W}_M} \, 
\widetilde{\mathfrak{f}}^{\rm tg}(v) := 
\mathfrak{f}^{\rm tg}(  \texttt{N}(\mathbf{a}) ).
\end{equation}
Then, we observe that
$$
D \widetilde{\mathfrak{f}}^{\rm tg}[v](h)
:=
\frac{1}{\epsilon}
\lim_{\epsilon\to 0} 
\frac{\widetilde{\mathfrak{f}}^{\rm tg}(v+\epsilon h) -
\widetilde{\mathfrak{f}}^{\rm tg}(v) }{\epsilon}
=
\nabla E(\mathbf{a}) \cdot \mathbf{h},
$$
for all $ v = \sum_{i=1}^M (\mathbf{a})_i \phi_i, h = \sum_{i=1}^M (\mathbf{h})_i \phi_i \in \mathcal{W}_M$.

In \cite{neuberger2009sobolev}, the author  considers a variational counterpart of the gradient-descent technique \eqref{eq:wrong_gradient_descent} for \eqref{eq:optimization_variational}:
 $$
 (
v^{(k+1)}, \phi)_H
=
 (
v^{(k)},  \phi)_H
- \gamma^{(k)}
\,
D \mathfrak{f}^{\rm tg}[v^{(k)}] 
(\phi),
\quad
k=1,2,\ldots,
\qquad
\forall \, \phi \in \mathcal{W}_M.
 $$
If we introduce the matrix 
 $\mathbf{A}\in \mathbb{R}^{M\times M}$ such that $(\mathbf{A})_{i,j} = (\phi_j,\phi_i)_H$ for $i,j=1,\ldots,M$, the latter 
 corresponds  in the algebraic update
\begin{equation}
\label{eq:good_gradient_descent_VFs}
\mathbf{a}^{(k+1)}
=
\mathbf{a}^{(k)}
- \gamma^{(k)}
\mathbf{A}^{-1}
\left( 
\nabla E(\mathbf{a}^{(k)})
\right),
\quad
k=1,2,\ldots.
\end{equation}

We notice that \eqref{eq:good_gradient_descent_VFs}
critically depends on the choice of the inner product $(\cdot, \cdot)_H$ but it
 is independent of the choice of the basis of $\mathcal{W}_M$.
In \cite{labatut2025non}, 
following \cite{beg2005computing},
Labatut considers 
$\mathfrak{f}_{\rm pen}(v) = \frac{1}{2} \|  \mathcal{K} v\|_{L^2(\Omega)}^2$ and $(\cdot, \cdot)_H  = ( \mathcal{K} \cdot , \mathcal{K} \cdot      )_{L^2(\Omega)}$ where $\mathcal{K}$ is a suitable differential operator that ensures 
$H\subset \mathcal{V}_0(\Omega)$
We also notice that  \eqref{eq:wrong_gradient_descent} 
and \eqref{eq:good_gradient_descent_VFs} 
coincide for $H$-orthonormal bases of $\mathcal{W}_M$.

\subsection{Bijectivity and approximation properties of vector flows}
\label{sec:VB_maps_analysis1}
We present a thorough mathematical analysis of the properties of VFs; the proofs exploit  well-known results in mathematical analysis and  calculus of variations, 
and are presented in section \ref{sec:vector_flows}.
Proposition \ref{th:bijectivity_flows} 
shows that for any $\mathbf{a}\in \mathbb{R}^M$ the VF
$\texttt{N}(\mathbf{a})$ \eqref{eq:VB_maps} is a diffeomorphism: the ansatz 
\eqref{eq:VB_maps}  hence satisfies \eqref{eq:desiderata_mapping_space} without the need for introducing a penalty function. 

\begin{proposition}
\label{th:bijectivity_flows}
Let $\Omega$ be a bounded Lipschitz domain and let 
$v \in \mathcal{V}_0(\Omega)$. Define $\Phi: \overline{\Omega} \to \mathbb{R}^n$ such that
$\Phi=F[v]$.
 Then, $\Phi$ is a diffeomorphism (one-to-one and onto) in $\overline{\Omega}$.
\end{proposition}
 
Proposition \ref{th:approx_flows}
shows that any element of ${\rm Diff}_0(\overline{\Omega})$ can be expressed as the flow of a  vector field in $\mathcal{V}_0(\Omega)$:
this implies that 
if the sequence $\phi_1,\ldots, \phi_M,\ldots$ forms a  basis of $\mathcal{V}_0$
the VF ansatz 
\eqref{eq:VB_maps} 
is dense in ${\rm Diff}_0(\overline{\Omega})$
 for $M\to \infty$ --- see Proposition \ref{th:approximation_flows} for a more precise statement.
Therefore, we find that the unregularized problem
$\min_{\mathbf{a}  \in \mathbb{R}^M}
 \mathfrak{f}_{\mu}^{\rm tg}(\texttt{N}( \mathbf{a} ))$ converges to 
 $\min_{ \Phi \in {\rm Diff}_0(\overline{\Omega})   }
 \mathfrak{f}_{\mu}^{\rm tg}(\Phi)$ for $M\to \infty$.

\begin{proposition}
\label{th:approx_flows}
 Let  $\Omega$ 
 be a Lipschitz domain and let $\Phi \in {\rm Diff}_0(\overline{\Omega})$. Then, there exists $v \in \mathcal{V}_0(\Omega)$ such that   $\Phi=F[v]$.
\end{proposition}

Next result shows that the Jacobian determinant of VFs  is strictly positive in $\Omega$.
We remark that Proposition \ref{th:evolution_jacobian} is a well-known result in differential geometry that is of major importance in continuum mechanics
(see, e.g.,  \cite[Proposition 5.4, Chapter 1]{marsden1983mathematical}).
Note that in 
\eqref{eq:evolution_jacobian} we use the symbol $\nabla_x$ to indicate derivatives of the velocity field $v=v(x,t)$.
Finally, we observe that Proposition \ref{th:evolution_jacobian} provides a constructive way of building elements of ${\rm Diff}(\Omega) \setminus {\rm Diff}_0(\Omega)$:
let $\Omega$ be the three-dimensional unit ball centered in the origin and consider $\Phi:\Omega \to \Omega$ such that $\Phi(\xi) = -\xi$. By inspection, we verify that $\Phi$ is a diffeomorphism; however, $J={\rm det}(\nabla \Phi) = -1<0$. We conclude that 
$\Phi \in {\rm Diff}(\overline{\Omega}) \setminus {\rm Diff}_0(\overline{\Omega})$. For two-dimensional domains, it suffices to consider $\Phi(\xi) = {\rm vec}(\xi_1,-\xi_2)$.

\begin{proposition}
\label{th:evolution_jacobian}
 Let  $\Omega$ 
 be a Lipschitz domain and let
 $v\in \mathcal{V}_0(\Omega)$.
Denote by $X$ the flow of the vector field $v$ \eqref{eq:flow_diffeomorphisms}.
  Then, the Jacobian determinant $J(\xi,t):={\rm det} \left( \nabla  X(\xi,t) \right)$ satisfies
  $ J(\xi,t=0) = 1$ and
  \begin{equation}
  \label{eq:evolution_jacobian}
  \frac{\partial}{\partial t} \log J(\xi,t) = \nabla_x \cdot v (X(\xi,t),t),
\end{equation}   
which implies
$$
J(\xi,t) = {\rm exp}
\left(
\int_0^t \nabla_x \cdot v(X(\xi,s),s) \, ds
\right),
\qquad
\xi \in \Omega, t\in [0,1].
$$
The latter implies that $J$ is strictly positive in $\Omega$ (that is, $\Phi$ is orientation-preserving).
\end{proposition}

We conclude this section by recalling  that the flow depends continuously on data.
Here, we rely on the Cauchy-Lipschitz theorem for Lipschitz velocity fields: the analysis can be extended to significantly less regular transport fields through the DiPerna-Lions theory \cite{diperna1989ordinary}.

\begin{proposition}
\label{th:approximation_flows}
Let $v,w\in \mathcal{V}_0(\Omega)$, and let $L$ be the Lipschitz constant of the velocity $v$. Then, we have
\begin{equation}
\label{eq:continuity_data}
\| F[v] - F[w]  \|_{L^\infty(\Omega)}
\leq 
\frac{e^L - 1}{L}
\| v -w  \|_{L^\infty(\Omega \times (0,1))}.
\end{equation}
\end{proposition}

Proposition \ref{th:approximation_flows} is of paramount importance for numerical approximation. Given 
$\Phi = F[v] \in {\rm Diff}_0(\overline{\Omega})$ 
and the approximation space 
$\mathcal{W}_M:= {\rm span} \{\phi_i \}_{i=1}^M \subset \mathcal{V}_0(\Omega)$,  we find that
\begin{equation}
\label{eq:approximation_result_explained}
\inf_{\mathbf{a}\in \mathbb{R}^M} 
\| \Phi - \texttt{N}(\mathbf{a})  \|_{L^\infty(\Omega)}
\leq
\frac{e^L - 1}{L}
\inf_{\phi \in \mathcal{W}_M} 
\| v - \phi \|_{L^\infty(\Omega\times (0,1))},
\end{equation}
where $v\in \mathcal{V}_0(\Omega)$ is the velocity field associated with $\Phi$ --- i.e., $\Phi=F[v]$ --- and   $L=L(v)>0$ is the corresponding Lipschitz constant.
Note that the multiplicative constant in \eqref{eq:approximation_result_explained} does not depend on the Lipschitz constants of the elements of $\mathcal{W}_M$.

\subsection{Computation of the gradient of the target function}
In view of the application of gradient-based methods for
 \eqref{eq:tractable_optimization_based_registration}, we 
  provide explicit expressions for the gradient of the $M$-dimensional function $E: \mathbf{a} \mapsto \mathfrak{f}^{\rm tg}(\texttt{N}(\mathbf{a}))$. 
We state upfront that our analysis   exploits the fact that $\mathbf{a}$ is finite-dimensional: we refer to \cite[Theorem 8.10]{younes2010shapes} and \cite{polzin2018large} for the extension to the infinite-dimensional setting.

\begin{proposition}
\label{th:derivative_computations_VB}
The derivative of the VF \eqref{eq:VB_maps} with respect to the coefficients $\mathbf{a}$ is given by:
\emph{
\begin{equation}
\label{eq:derivative_VB}
\frac{\partial  \texttt{N} }   {\partial a_i}
(\xi; \mathbf{a})
=
\nabla X(\xi, 1; \mathbf{a}) 
\int_0^1 \, 
\left(
\nabla X(\xi, \tau; \mathbf{a}) 
\right)^{-1}
\phi_i \left( X(\xi, \tau; \mathbf{a}), \tau \right) \, d\tau,
\end{equation}
}
for $i=1,\ldots,M.$ 
The  derivative of $E$ is then obtained by plugging 
\eqref{eq:derivative_VB} in
 \eqref{eq:sensitivity_method_b}.\end{proposition}

Computation of \eqref{eq:sensitivity_method_b} requires the solution to the nonlinear ODE system \eqref{eq:flow_diffeomorphisms} and to $n^2$ linear ODE systems for the calculation of $\nabla X$ 
(see Eq. \eqref{eq:ODE_gradient} for the explicit formula). 
An alternative approach, which has been extensively used in PDE-constrained optimization
(see,  e.g., \cite{jameson1988aerodynamic,pironneau2005optimal}) and also image registration \cite{polzin2018large}, relies on the solution to an adjoint equation. This second approach   only requires the solution to $n$ additional linear ODE systems.
In the next result, we provide the expression of the adjoint method for the two target functions introduced above:
notice that the computation of the gradient using the adjoint method requires the solution to 
a nonlinear  forward problem for the computation of $X$ (cf. \eqref{eq:flow_diffeomorphisms}), the solution to $n$ linear backward problems 
(cf.  \eqref{eq:adjoint_method_distributed_b} and 
\eqref{eq:adjoint_method_pointwise_b})
 for the computation of $\Lambda$ and a simple function  evaluation  (cf.  \eqref{eq:adjoint_method_distributed_a} or
\eqref{eq:adjoint_method_pointwise_a}).

\begin{proposition}
\label{th:adjoint}
Consider 
$\mathfrak{f}^{\rm tg}= \mathfrak{f}_{(1)}^{\rm tg}$
in \eqref{eq:distributed_sensor} and set
\emph{$\Psi(\cdot; \mathbf{a}) := \big( u \circ \texttt{N}( \mathbf{a}) 
$
$
-
\zeta_\Phi( \mathbf{a}) \big)
\nabla_x u \circ \texttt{N}( \mathbf{a})$}
(cf. \eqref{eq:sensors_der}). Then,
\begin{subequations}
\label{eq:adjoint_method_distributed}
\begin{equation}
\label{eq:adjoint_method_distributed_a}
\frac{\partial E}{\partial a_k} (\mathbf{a})= 
\int_{\Omega \times (0,1) }
\Lambda(\xi,t; \mathbf{a})
\cdot 
\phi_k\left(   X(\xi,t; \mathbf{a}) , t \right) \, d\xi dt,
\quad
k=1,\ldots,M,
\end{equation}
where the adjoint $\Lambda(\mathbf{a}): \Omega \times [0,1] \to \mathbb{R}^n$ satisfies:
\begin{equation}
\label{eq:adjoint_method_distributed_b}
\left\{
\begin{array}{ll}
\displaystyle{
\frac{\partial \Lambda}{\partial t}(\xi,t; \mathbf{a}) = - \nabla_x v \left( X(\xi,t),t ; \mathbf{a} \right) ^\top \Lambda(\xi,t; \mathbf{a}) 
}
& t\in (0,1),
\\[3mm]
\Lambda(\xi, 1; \mathbf{a}) = 
\Psi(\xi; \mathbf{a})
&
\end{array}
\right.
\qquad
\forall \, \xi \in \Omega.
\end{equation}
\end{subequations}

Consider 
$\mathfrak{f}^{\rm tg}= \mathfrak{f}_{(2)}^{\rm tg}$.  Then, we find
\begin{subequations}
\label{eq:adjoint_method_pointwise}
\begin{equation}
\label{eq:adjoint_method_pointwise_a}
\frac{\partial E}{\partial a_k} (\mathbf{a}) \, = \, 
\sum_{i=1}^{N_0}
\int_0^1
\Lambda_i(t; \mathbf{a}) \, \cdot \,
\phi_k \left(   X(\xi_i, t; \mathbf{a})    \right)  \, dt,
\quad
k=1,\ldots,M;
\end{equation}
where the time-dependent functions  $\Lambda_1(\mathbf{a}), \ldots,\Lambda_{N_0}(\mathbf{a})  :  [0,1] \to \mathbb{R}^n$ satisfy :
\emph{
\begin{equation}
\label{eq:adjoint_method_pointwise_b}
\left\{
\begin{array}{ll}
\displaystyle{
\frac{d \Lambda_i}{d t}(t; \mathbf{a}) 
\, = \, 
 - \nabla_x v \left( X(\xi_i,t),t ; \mathbf{a} \right) ^\top \, \Lambda_i(t; \mathbf{a})
}
& t\in (0,1),
\\[3mm]
\displaystyle{
\Lambda_i( 1; \mathbf{a}) 
\, = \, 
\texttt{N}(\xi_i; \mathbf{a}) - 
\sum_{j=1}^{N_1}
P_{i,j}  y_j
}
\\
\end{array}
\right.
i=1,\ldots,N_0.
\end{equation}
}
\end{subequations}
\end{proposition}

\begin{remark}
The expressions in Proposition \ref{th:adjoint}
rely on the adjoint of the continuous equations:
in practice, the solution $X$ to 
\eqref{eq:flow_diffeomorphisms},
the gradient $\nabla X$ that enters in
\eqref{eq:derivative_VB}, and the adjoint states 
\eqref{eq:adjoint_method_distributed_b}-\eqref{eq:adjoint_method_pointwise_b}
 are obtained through standard Runge Kutta methods that are not dual consistent. Therefore, the estimated gradient 
 is only consistent in the limit.
\end{remark}

\subsection{Bijectivity and approximation properties of compositional maps in polytopes}
\label{sec:DB_maps_analysis}

We first introduce a class of domains that are important for the analysis. 
This class is the higher-dimensional analogue of the family of polygonal domains introduced in \cite{taddei2025compositional}. They consist in bounded Lipschitz domains that  admit a ``simplicial'' mesh. 
{Below, $\card(A)$ denotes the cardinality of the set $A$}.

\begin{Definition}[simplex]
A {closed} subset $\sigma$ of $\mathbb{R}^n$ is called a $k$-dimensional simplex (or $k$-simplex) if it can be written as the convex hull of a finite set $\{x_1,\dots,x_{k+1}\}$ and if it is not contained in a $k-1$-dimensional affine subspace. 
For any non-empty $I\subset \{x_1,\dots,x_k\}$, the convex hull of $I$ is called a face of $\sigma$. It is  a $(\card(I)-1)$-simplex.
\end{Definition}

We observe that a $0$-simplex is a point in $\mathbb{R}^n$, a $1$-simplex a segment, a $2$-simplex a triangle, 
a $3$-simplex a tetrahedron;
{By construction, any face of a simplex belongs to the boundary   $\partial \sigma$. Below, we provide the definition of polyhedra and polyhedral domains (see, e.g., \cite{ciarlet2002finite}).}

\begin{Definition}[polyhedron]
\label{def:polyhedron}
A subset  $K$ of $\mathbb{R}^n$ is a bounded polyhedron if there exists a finite set of $n$-simplices of $\mathbb{R}^n$ 
$\sigma_1,\dots, \sigma_\ell$ such that 
\begin{enumerate}
\item 
the union of the simplices $\sigma_i$ is $K$;
\item 
for any $1\leq i< j\leq \ell$,  $\sigma_i\cap\sigma_j$ is either empty or it is an entire face of both of them.
\end{enumerate}
{The family 
$\mathcal{T} :=\{ \sigma_1,\dots,\sigma_\ell \}$ is called a triangulation of $K$;
we further denote by $\mathcal{F}_k$ the family of $k$-dimensional faces of $\sigma_1,\dots,\sigma_\ell$, for $k=0,1,\ldots,n-1$.}
\end{Definition}

\begin{Definition}[polyhedral domain]
\label{def:polyhedral_domain}
The open set 
 $\Omega\subset \mathbb{R}^n$ is a bounded polyhedral domain  if $\Omega$ is a connected Lipschitz domain of $\mathbb{R}^n$ and $\oo$ is a bounded polyhedron.
\end{Definition}

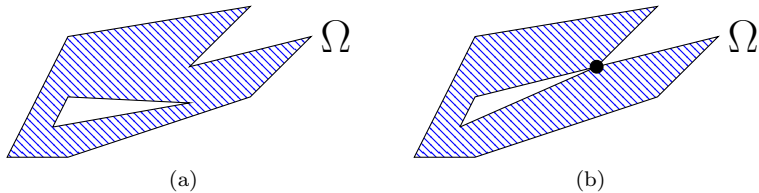
\begin{figure}[h!]
\centering
\subfloat[]{
\centering
\begin{tikzpicture}[scale=0.8]
\linethickness{0.3 mm}
\linethickness{0.3 mm}

\filldraw[pattern=north west lines, pattern color=blue] 
(0,0)--(3,1)--(4,2)--(2,1.5)--(3,2.5)--(0,2)--(-1,0)--(0,0);

\filldraw[fill=white]  (-0.25,0.5)--(2,0.9)--(0,1)--(-0.25,0.5);



\coordinate [label={right:  {\LARGE {$\Omega$}}}] (E) at (4, 2) ;
\end{tikzpicture}
}
~~~
\subfloat[]{
\begin{tikzpicture}[scale=0.8]
\linethickness{0.3 mm}
\linethickness{0.3 mm}

\filldraw[pattern=north west lines, pattern color=blue] 
(0,0)--(3,1)--(4,2)--(2,1.5)--(3,2.5)--(0,2)--(-1,0)--(0,0);

\filldraw[fill=white]  (-0.25,0.5)--(2,1.5)--(0,1)--(-0.25,0.5);

\filldraw (2,1.5) circle (3pt);


\coordinate [label={right:  {\LARGE {$\Omega$}}}] (E) at (4, 2) ;
\end{tikzpicture}
}
\caption{interpretation of Definitions
\ref{def:polyhedron} and \ref{def:polyhedral_domain}.
(a) example of polyhedral domain.
(b) example of polyhedron that is not Lipschitz.
}
\label{fig:geometric_definition}
\end{figure}

{Polyhedra are unions of a finite number of non-overlapping simplices. We observe that the boundary of a polyhedron is not guaranteed to be Lipschitz (cf. 
Figure \ref{fig:geometric_definition}(b) from \cite{taddei2025compositional}).}

Proposition \ref{th:bijectivity} shows that the CM \eqref{eq:DB_ansatz_polytope} is a bijection in $\Omega$, provided that 
$\mathfrak{f}_{\rm pen}(\mathbf{a}) = - \inf_{x\in \Omega} {\rm det} (\nabla \texttt{N}(\mathbf{a})) < 0$.
On the other hand, Proposition \ref{th:approximation_DB} shows that CMs are dense in ${\rm Diff}_0(\overline{\Omega})$, provided that the basis $\varphi_1,\varphi_2,\ldots, $ is dense in $\mathcal{U}_0(\Omega)$.
Proposition \ref{th:approximation_DB}  can also be used to show that ${\rm Diff}_0(\overline{\Omega})$ is strictly contained in ${\rm Diff}(\overline{\Omega})$. Consider $\Omega=(-1,1)^3$ and consider $\Phi(\xi) = {\rm vec}(\xi_1,-\xi_2,-\xi_3)$; clearly, $\Phi$ is a diffeomorphism with positive Jacobian determinant,  but $\varphi:=\Phi - \texttt{id}$ does not  belong to $ \mathcal{U}_0(\Omega)$; therefore, $\Phi \in {\rm Diff}(\overline{\Omega}) \setminus {\rm Diff}_0(\overline{\Omega})$.

\begin{proposition}
\label{th:bijectivity}
 Let  $\Omega$ 
 be 
a  polyhedral domain of $\mathbb{R}^n$.
Consider the vector-valued function
\emph{$\Phi = \texttt{id} + \varphi$} with 
$\varphi \in \mathcal{U}_0(\Omega)$. Then $\Phi$ is a bijection (one-to-one and onto) in $\Omega$ if $\min_{x \in \overline{\Omega}} 
J(\Phi ) > 0$.
\end{proposition}

\begin{proposition}
\label{th:approximation_DB}
Let 
$\Omega$ be a 
 polyhedral domain of $\mathbb{R}^n$
and let 
$\Phi \in {\rm Diff}_0(\overline{\Omega})$. Then, 
\emph{$\Phi = \texttt{id}+\varphi$} with $\varphi \in \mathcal{U}_0(\Omega)$.
\end{proposition}

\begin{remark}
[Extension to curved domains.]
Since the composition of Lipschitz bijections is a Lipschitz bijection, we can also exploit Proposition \ref{th:bijectivity} to show that 
{$\texttt{N}(\mathbf{a})$} in \eqref{eq:DB_ansatz} is a bijection in $\Omega$ if 
{$J(\texttt{N}_{\rm p}(\mathbf{a}))$} is strictly positive in $\overline{\Omega}_{\rm p}$.
As discussed in \cite{taddei2025compositional} for $n=2$-dimensional domains, the nonlinear generalization \eqref{eq:DB_ansatz} of \eqref{eq:DB_ansatz_polytope} for curved domains is not dense in $ {\rm Diff}_0(\overline{\Omega})$. The same negative result readily applies to three-dimensional domains. In  \cite{taddei2025compositional}, we show that it is possible to enhance the expressivity of CMs by considering compositions of CMs associated with different polytopes (cf. Lemma 3.3 in \cite{taddei2025compositional}).  In this work, we do not discuss this issue any further.
\end{remark}

\subsection{Discussion}
Evaluation of VFs involves the solution to a system of ordinary differential equations (ODEs) for each point, while CMs simply require three function evaluations --- furthermore, the evaluation of $\Psi^{-1}$ can be pre-computed where needed before 
solving \eqref{eq:tractable_optimization_based_registration}
as it is independent of the mapping coefficients $\mathbf{a}$. VFs are hence much more expensive to evaluate.

Direct computation of 
$\frac{\partial \texttt{N}}{\partial \mathbf{a}}  (\mathbf{a})$ for VFs at the discrete level is computationally expensive for typical (e.g., Runge Kutta, linear multi-step) time-integration schemes:
on the other hand, the evaluation of   the \emph{continuous} gradient of the objective function
(cf. Proposition \ref{th:derivative_computations_VB})
 can be performed with  roughly the same computational cost as the cost of evaluating $\texttt{N} (\mathbf{a})$.

The VF \eqref{eq:VB_maps}
is guaranteed to be  a diffeomorphism of $\overline{\Omega}$ (cf. Proposition \ref{th:bijectivity_flows}) for all $\mathbf{a}\in \mathbb{R}^M$; furthermore, any element of ${\rm Diff}_0(\overline{\Omega})$  can be expressed in the form \eqref{eq:flow_diffeomorphisms} (cf. Proposition \ref{th:approx_flows}),
which implies that  
VFs \eqref{eq:VB_maps} are dense in  ${\rm Diff}_0(\overline{\Omega})$ for $M\to \infty$ --- provided that the sequence $\phi_1,\phi_2,\ldots$ is dense in $\mathcal{V}_0$ (cf. \eqref{eq:approximation_result_explained}). Note, however, that approximate integration of the ODE system in \eqref{eq:flow_diffeomorphisms} does not preserve the bijectivity constraint: it is hence necessary to accurately integrate \eqref{eq:flow_diffeomorphisms}, particularly in the proximity of the boundary of $\Omega$.

CMs of the form \eqref{eq:DB_ansatz_polytope}-\eqref{eq:DB_ansatz} are bijective if   the minimum value of the Jacobian determinant  is strictly positive
(cf. Proposition \ref{th:bijectivity}); therefore,  the penalty function in \eqref{eq:tractable_optimization_based_registration} for CMs  is highly non-linear.
For parametric registration tasks,  it is hence difficult to ensure 
bijectivity 
 for out-of-sample parameters with 
standard regression algorithms (see
\cite{taddei2025compositional}).
Furthermore,  CMs are dense in ${\rm Diff}_0(\overline{\Omega})$ if and only if $\Omega$ is a polygonal domain
(cf. Proposition \ref{th:approximation_DB}): this implies that for curved domains their expressive power critically depends on the choice of the polygonal domain $\Omega_{\rm p}$ that might be difficult to select in practice \cite{taddei2025compositional}.

In conclusion, our analysis reveals a fundamental trade-off: while VFs constitute the only mathematically rigorous choice 
for generic curved domains 
--- as they are dense in ${\rm Diff}_0(\overline{\Omega})$ and allow a simple enforcement of the  bijectivity constraint --- this theoretical 
consistency comes at the price of the significantly higher computational cost associated with flow integration. 
Conversely, CMs offer a computationally attractive algebraic alternative, but their density in 
${\rm Diff}_0(\overline{\Omega})$
is strictly limited to polyhedral domains, as they necessitate an auxiliary transformation for curved geometries.

\section{Choice of the truncated basis}
\label{sec:truncated_basis}
The effectiveness of VFs and CMs strongly depends on the choice of the expansions in \eqref{eq:VB_maps} and \eqref{eq:DB_ansatz}. In the previous works
\cite{labatut2025non,taddei2025compositional}, we relied on traditional finite element (FE) discretizations:
this choice leads to very high-dimensional expansions and thus to very expensive computations. To address this issue,  we propose two distinct strategies for the construction   of empirical expansions, which reduce the dimensionality of the optimization problem \eqref{eq:tractable_optimization_based_registration} \emph{a priori}, and, in perspective, enable the use of  more sophisticated second-order optimization routines that are particularly effective for minimization problems of moderate size. 

\subsection{Modal expansions for registration}
To simplify the presentation and the  mathematical analysis, given the Hilbert space $H \subset \mathcal{U}_0(\Omega)$ (or $H \subset \mathcal{V}_0(\Omega)$ for VFs), we introduce the high (but finite) dimensional space $H_{\rm hf}:=  {\rm span}\{ \varphi_i  \}_{i=1}^N\subset H$. Given $u\in H_{\rm hf}$, we denote by $\mathbf{u}\in \mathbb{R}^N$ the corresponding vector of coefficients $u= \sum_{i=1}^N (\mathbf{u})_i \varphi_i$. Next, we introduce the high-dimensional minimization problem
\begin{equation}
\label{eq:target_optimization}
\mathfrak{E}:= \min_{\mathbf{u}\in \mathbb{R}^N} f(\mathbf{u}) + \xi \| \mathbf{u} \|_{\mathbf{A}}^2
\end{equation}
where $\xi>0$, $ \| \mathbf{u} \|_{\mathbf{A}}^2 = \mathbf{u}^\top \mathbf{A} \mathbf{u}$, $\mathbf{A}\in \mathbb{R}^{N\times N}$ is a symmetric positive semi-definite matrix, $f:\mathbb{R}^N \to \mathbb{R}$ is a Lipschitz function with Lipschitz constant
\begin{equation}
\label{eq:lipschitz_constant}
L(r):= \max_{ \mathbf{u},\mathbf{u}' \in \mathbb{R}^N , \mathbf{u} \neq \mathbf{u}',
\| \mathbf{u}  \|_{\mathbf{M}},
\| \mathbf{u}'  \|_{\mathbf{M}} \leq r
} 
\frac{|f(\mathbf{u})  -f(\mathbf{u}')|}{  \| \mathbf{u} \|_{\mathbf{M}}     },
\end{equation}
and $\mathbf{M}\in \mathbb{R}^{N\times N}$ is a symmetric positive definite matrix.
Note that \eqref{eq:target_optimization} reads as a minimization problem with quadratic Tykhonov regularization: we remark that the approaches in 
\cite{labatut2025non,taddei2025compositional} and also \cite{younes2010shapes} can be recast as in \eqref{eq:target_optimization}. 
Exploiting this notation, we aim to find a reduced-order basis (ROB) $\mathbf{W} \in \mathbb{R}^{N\times m}$ with $m\ll N$ such that
\begin{equation}
\label{eq:target_optimization_reduced}
\mathfrak{E}_m:= \min_{\mathbf{a}\in \mathbb{R}^m} f(\mathbf{W} \mathbf{a}) + \xi \| \mathbf{W} \mathbf{a} \|_{\mathbf{A}}^2
\approx \mathfrak{E}.
\end{equation}

Our first strategy consists in taking the eigenmodes of the generalized eigenproblem $\mathbf{A} \boldsymbol{\phi}  = \lambda  \mathbf{M} \boldsymbol{\phi}$ associated with the $m$ smallest eigenvalues, that is
\begin{equation}
\label{eq:empirical_eigen}
\mathbf{W} = \left[ 
\boldsymbol{\phi}_1,\ldots, \boldsymbol{\phi}_m
\right],
\quad
\mathbf{A} \boldsymbol{\phi}_i  = \lambda_i  \mathbf{M} \boldsymbol{\phi}_i,
\quad
i=1,\ldots,N; \;\;
0\leq \lambda_1 \leq \ldots \leq \lambda_N.   
\end{equation}
We remark that a similar strategy was proposed in
\cite[section 6]{myronenko2010point} in the context of point set registration in unbounded domains.

Next Lemma justifies the definition of $\mathbf{W}$ in \eqref{eq:empirical_eigen}:
the proof is straightforward and is provided in 
section \ref{sec:easy_results}.
The result suggests that if the eigenvalues of $\mathbf{M}^{-1} \mathbf{A}$ increase sufficiently fast we can consider a modal basis with $m \ll N$ at the price of solving a large-scale symmetric  eigenproblem; clearly, the cost of the latter can be amortized for many query (parametric) registration tasks. In the numerical experiments, we investigate the performance of this choice.
Notice that the bounds of Lemma \ref{lemma:eigen_justification} depend on the norm $r$ of the minimizer $\mathbf{u}^\star$: we envision that $r$ can be bound from above by introducing additional assumptions on the function $f$. We omit the details.

\begin{lemma}
\label{lemma:eigen_justification}
Let $\mathbf{u}^\star$ be a global minimizer of \eqref{eq:target_optimization} and denote by $P_W: \mathbb{R}^N \to {\rm col}(\mathbf{W})$ the $\|\cdot\|_{\mathbf{M}}$-orthogonal projection operator and set 
 $P_W^\perp = \mathbbm{1} - P_W$. Then, we have
 \begin{equation}
\label{eq:eigen_justification_1}
 \| P_W^\perp \mathbf{u}^\star \|_{\mathbf{M}}
\leq
\frac{L^\star}{\xi \lambda_{m+1}}, 
\quad
{\rm with} \;\;
L^\star = L(  \| \mathbf{u}^\star  \|_{\mathbf{M}} );  
 \end{equation}
and
 \begin{equation}
\label{eq:eigen_justification_2}
0\leq \mathfrak{E} - \mathfrak{E}_m \leq 
\dfrac{(L^\star)^2}{\xi \lambda_{m+1}}.
 \end{equation}
\end{lemma}

We also investigate an alternative strategy based on the solution to suitable linear problems with polynomial source terms. Let $a:H \times H \to \mathbb{R}$  be a suitable linear form and let $\mathbb{P}_p(\mathbb{R}^n)$ be the space of polynomials of $n$ variables of total degree less or equal to $p$; then, we define the generalized finite element space:
\begin{equation}
\label{eq:gfem_space}
\mathcal{W}_p := 
\left\{
u\in H_{\rm hf} \, : \,
a(u,v) = \int_{\Omega} f \cdot v \, dx + \int_{\partial \Omega} g\cdot v \, dx,
\;\;
\forall \, v \in H_{\rm hf}, \;\;
f,g \in [\mathbb{P}_p(\mathbb{R}^n)]^n
\right\}
\end{equation}
In the numerical experiments, we consider the form $a(u,v) =\mathbf{v}^\top \mathbf{A} \mathbf{u}$.
For VFs, 
we can simply consider the space-time tensorized extension of \eqref{eq:gfem_space}:
\begin{equation}
\label{eq:gfem_space_time}
\mathcal{W}_p^{\rm st} := 
\left\{
u(x,t) = u'(x) \ell(t) \, : \,
u' \in \mathcal{W}_p, \;\;
\ell\in \mathbb{P}_p(\mathbb{R})
\right\}.
\end{equation}

We conclude this section by briefly commenting on the imposition of the constraint
$u\cdot \mathbf{n}\big|_{\partial \Omega} = 0$. In \cite{taddei2025compositional}, the boundary constraint is enforced directly in the space, while in \cite{labatut2025non} the constraint is imposed through the Nitsche's method. For curved domains, the quality of the approximation of the normal has a significant impact on accuracy: in this work, we do not further address this issue.

\subsection{Numerical results: choice of the modal expansion for CMs}
\label{sec:numerics}
\begin{figure}[h!]
\centering
\subfloat[]{ 
\includegraphics[width=.4\textwidth]{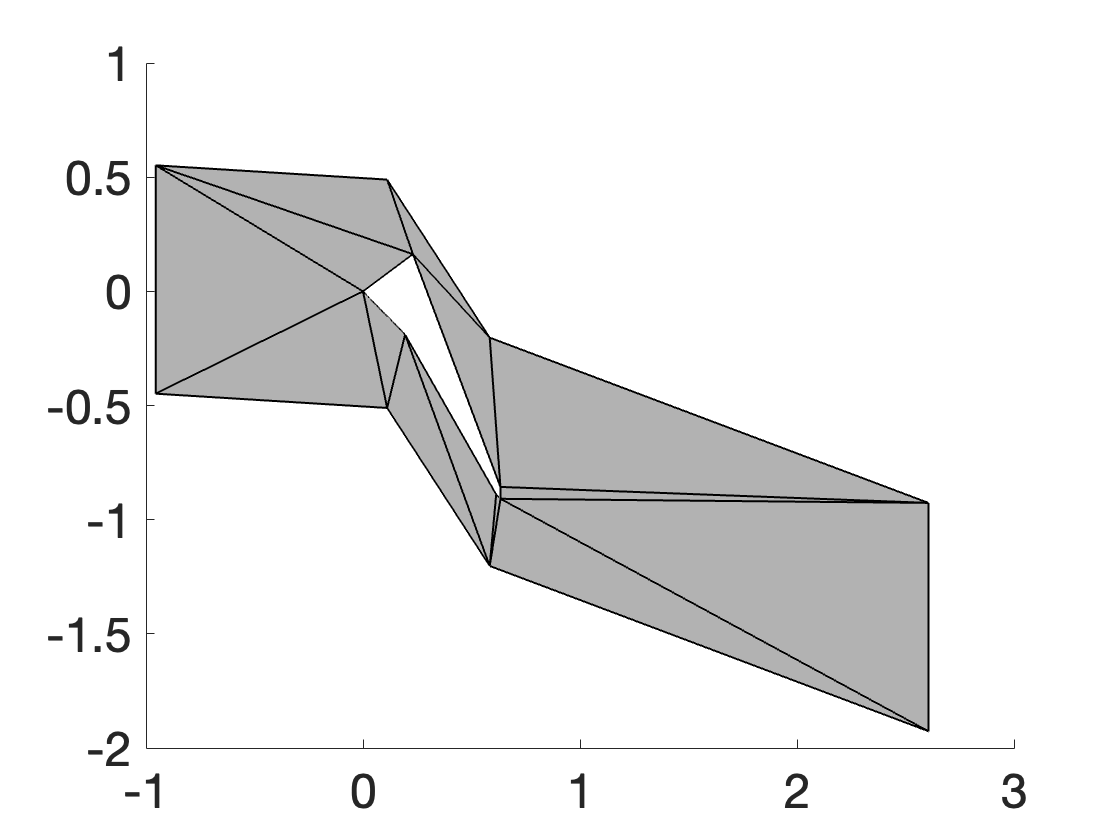}
}
~~
\subfloat[]{
\includegraphics[width=.4\textwidth]{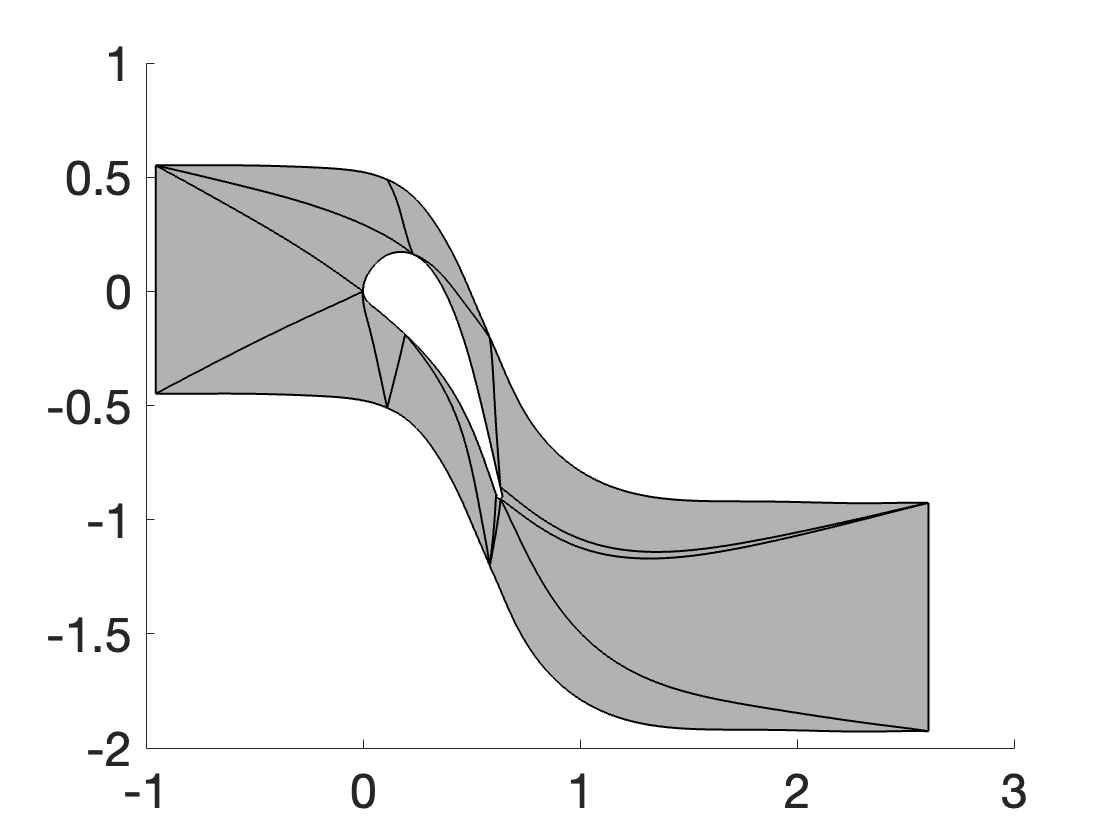}
}
\caption{Inviscid flow past an array of LS89 turbine blades.
(a) polytope $\Omega_{\rm p}$ and FE mesh;
(b) curved domain $\Omega$.
}
\label{fig:LS89_meshes}
\end{figure}

We consider the test case described in \cite[section 5]{taddei2025compositional}, which exploits compositional maps for registration. We refer to the original paper for a detailed presentation of the test case.
Figure \ref{fig:LS89_meshes} shows the polytope $\Omega_{\rm p}$ (cf. \eqref{eq:DB_ansatz}) and the FE mesh associated with the space $H_{\rm hf}$ ---
in the numerical experiments, we consider polynomials of degree $7$ ($N=784$) and we apply the greedy procedure  in \cite[Algorithm 4.1]{taddei2025compositional} based on a training set with $n_{\rm train}=66$ parameters to generate the FE solutions $\{ \varphi_\mu: \mu \in \mathcal{P}_{\rm train}\}$.

\begin{figure}[h!]
\centering
\subfloat[]{ 
\includegraphics[width=.4\textwidth]{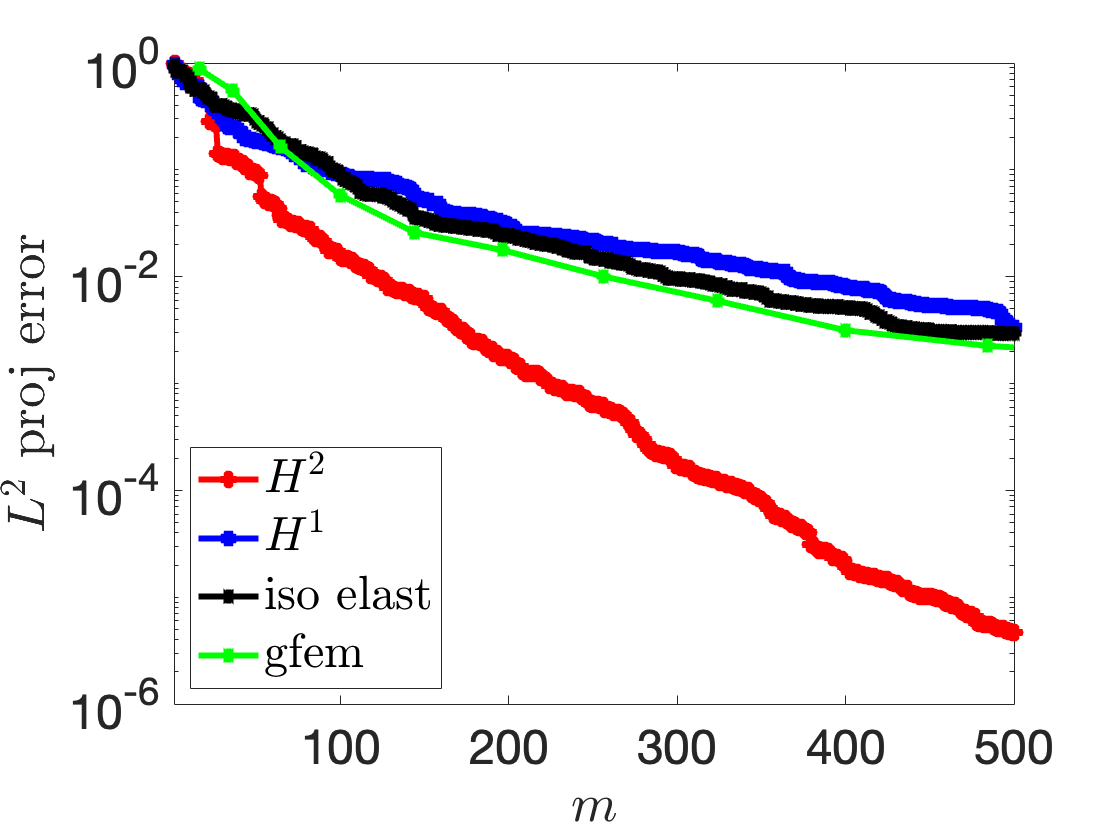}
}
~~
\subfloat[]{
\includegraphics[width=.4\textwidth]{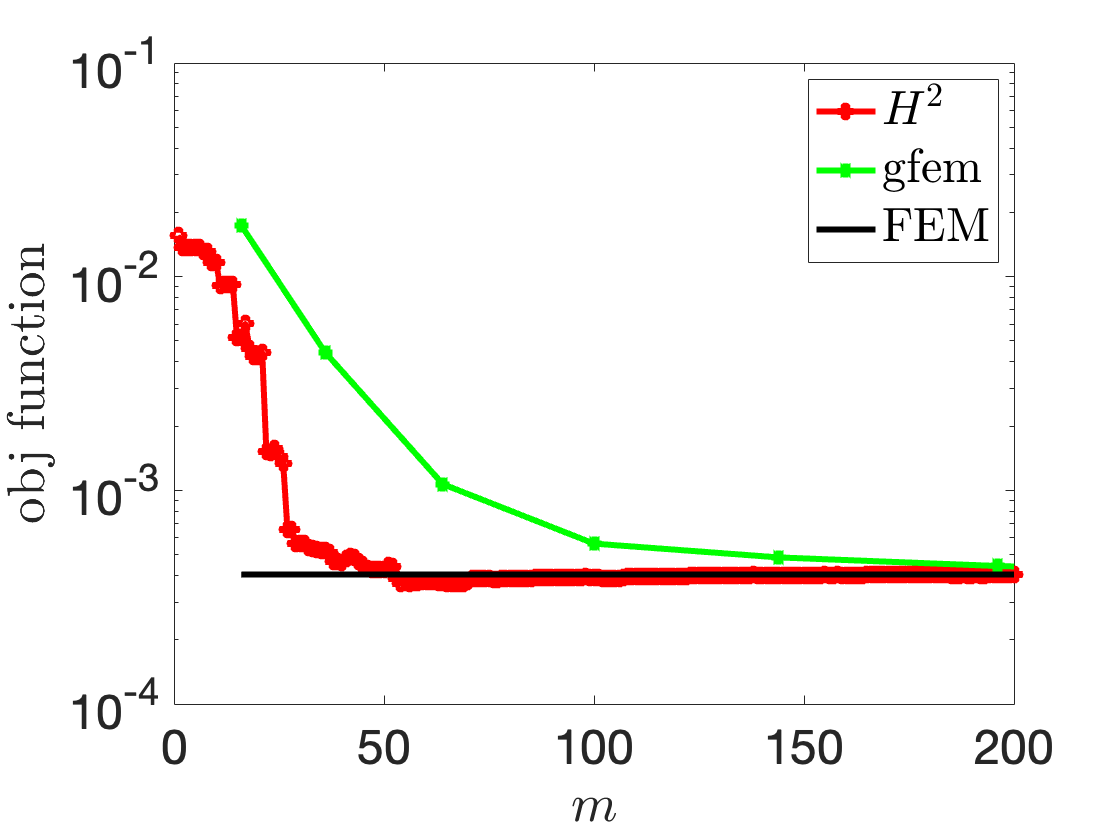}
}
\caption{Inviscid flow past an array of LS89 turbine blades.
(a) relative projection error $E_m^{\rm proj}$.
(b) worst-case  objective $E_m^{\rm obj}$ and comparison with FE solution.
}
\label{fig:LS89_results}
\end{figure}

Figure \ref{fig:LS89_results} shows the results. 
Figure \ref{fig:LS89_results}(a) shows the behavior of the worst-case $L^2$ projection error for the empirical bases of section \ref{sec:numerics}.
\begin{equation}
\label{eq:E_proj_ls89}
E_m^{\rm proj} = \max_{\mu\in \mathcal{P}_{\rm train}}
\dfrac{\|\varphi_\mu - P_{\mathcal{W}_m} \varphi_\mu  \|_{L^2(\Omega_{\rm p})}}{\|\varphi_\mu  \|_{L^2(\Omega_{\rm p})}};
\end{equation}
Figure \ref{fig:LS89_results}(b) shows the behavior of the objective function  
\begin{equation}
\label{eq:E_obj_ls89}
E_m^{\rm obj} = \max_{\mu\in \mathcal{P}_{\rm train}}
\mathfrak{f}_\mu^{\rm tg} ( P_{\mathcal{W}_m} \varphi_\mu ).
\end{equation}
We notice that the eigenspace associated with the $H^2$ seminorm leads to superior performance  if compared to the eigenspaces associated with the standard $H^1$
norm 
$$
\int_{\Omega_{\rm p}} \nabla \phi_i : \nabla v \, dx
=
\lambda_i 
\int_{\Omega_{\rm p}}   \phi_i \cdot v \, dx
\qquad
\forall \, v\in H_{\rm hf},
$$
and the  isotropic linear elasticity operator (with plane-strain assumption)
$$
\int_{\Omega_{\rm p}}
\frac{E\nu}{(1+\nu)(1-2\nu)}
\nabla_{\rm s} \phi_i : \nabla_{\rm s}  v \,
+
\frac{E}{2(1+\nu)}
(\nabla \cdot \phi_i) 
(\nabla \cdot v) 
dx
=
\lambda_i 
\int_{\Omega_{\rm p}}   \phi_i \cdot v \, dx
\;\;
\forall \, v\in H_{\rm hf},
$$
with $E=1$ and $\nu=\frac{1}{3}$.
We further notice that the performance associated with \eqref{eq:gfem_space}
(label ``gfem'')
are inferior to the performance of the $H^2$ eigenspace.
Finally, we remark that both the $H^2$ eigenspace method and the gfem method achieve comparable performance in terms of the target function for modest values of $m$.

 \section{Proofs of Propositions \ref{th:bijectivity} and \ref{th:approximation_DB}}
  \label{sec:proofs}

We rely on several topological results that do not exploit the Lipschitz regularity of the domain nor the differentiability of the mapping $\Phi$: for completeness, we report some definitions and results that are used in the proofs.
Given the set $U \subset \mathbb{R}^n$, we denote by $\partial U$ the boundary of $U$, which is the  closure of $U$ minus the interior of $U$; we further denote by 
 $\mathcal{O}_{U}$ the induced (or subspace) topology on $U$,
$
\mathcal{O}_{U} := \{A \cap U: \, A \, {\rm is \, open \, in} \;  \mathbb{R}^d \}.
$
It is possible to show that the ordered pair
$\mathcal{T}_{U} = (U, \mathcal{O}_{U})$ is a topological space. We say that $B$ is open in $U$ if $B \in \mathcal{O}_{U}$; similarly, we say that $B$ is closed in $U$ if the complement of $B$ in $U$,
$B^{\rm c}:= U \setminus B$,
belongs to  $\mathcal{O}_{U}$.
We further say that $\mathcal{T}_{U}$ is connected if it cannot be represented as the union of two disjoint non-empty open subsets; we say that $\mathcal{T}_{U}$ is path-connected if there exists a path joining any two points in $U$; finally, we say that 
$\mathcal{T}_{U}$ is simply-connected if $\mathcal{T}_{U}$ is path-connected and every path between two points can be continuously transformed into any other such path while preserving the endpoints.  
Given the topological spaces $X,Y$, $f: X\to Y$ is a homeomorphism if (i) it is a bijection (one-to-one and onto), (ii) it is continuous, and (iii) the inverse $f^{-1}$ is continuous (i.e., $f$ is an open mapping).

Theorems \ref{th:custom_theorem} and \ref{th:connected_sets} are standard results in topology that can be found in \cite[Chapter 3]{munkres2013topology} (cf. \S 23 and Exercise 6);
 
\begin{theorem}
\label{th:custom_theorem}
In a topological space $T$, any connected subset of $T$ that meets
both a set $C$ and its complement $T \setminus C$ intersects the boundary of $C$.
\end{theorem} 


\begin{theorem}
\label{th:connected_sets}
A topological space $\mathcal{T}_{U}$ is connected if and only if the only open and closed sets are the empty set and $U$.  
\end{theorem}

Theorem \ref{th:brouwer} is Brouwer's invariance of domain theorem and  can  be found in \cite[Chapter 10]{munkres2013topology}. Its proof is much more involved than the latter ones. 
\begin{theorem}
\label{th:brouwer}
If $U$ is an open subset of $\mathbb R^n$ and $f:U\rightarrow \mathbb {R}^{n}$ is an injective continuous map, then $V := f( U )$ is open in $\mathbb {R}^{n}$ and $f$ is a homeomorphism between $U$  and $V$.
\end{theorem}

A subset $U$ of $\mathbb{R}^n$ is a domain with continuous boundary if $U$ is open in $\mathbb{R}^n$ and if its boundary $\partial U$ is locally the graph of a continuous function.
{Note that  the boundary of domains of this class might feature cusps; however, they cannot feature self-intersecting boundaries (e.g., Figure \ref{fig:geometric_definition}(b)).}
If $U\subset \mathbb{R}^n$ is a domain with continuous boundary, then $\partial U$ is \emph{locally flat}  \cite{Brown}  i.e. for any point $\xi$ in $\partial U$ there exists an open set $N\ni \xi$  and a homeomorphism $h:N\rightarrow h(N)\subseteq \mathbb{R}^n$ such that
{$h(N\cap \partial U)\subseteq  \mathbb{R}^{n-1}
\times (-\infty, 0) $} and
 $h(N\cap \partial U)\subseteq 
\mathbb{R}^{n-1} \times  \{0\}$. 
Locally-flat domains have an important property that was proved by   Brown (cf. 
\cite[Theorem 3]{Brown}).

\begin{theorem}
\label{th:brown}
Let $\Omega$ be an open set with continuous (and thus locally-flat) boundary.
For any connected component 
$\Sigma$
of $\partial \Omega$, there exists a homeomorphism $h_\sigma$ from $\Sigma \times (-1,1)$ to a neighborhood of $\Sigma$ in $\mathbb{R}^n$ such that for any $\xi\in \Sigma$, $h_\sigma(\xi,0)=\xi$. Such a map  is called a tubular neighborhood of $\Sigma$.
\end{theorem}

Theorem \ref{th:hadamard2} is a generalization of Hadamard's global inverse function theorem and is proven in  \cite{Ho}.  As pointed out by the author of  \cite{Ho}, this version can be used on manifolds with boundary and in particular on bounded domains with Lipschitz boundary.
We recall that an Hausdorff topological space is a topological space where distinct points have distinct neighborhoods; furthermore, the function $f:X \rightarrow  Y$ is   proper if the preimage of any compact set in $Y$ is compact in $X$.

\begin{theorem}
\label{th:hadamard2}
Let $X$ and $Y$ be Hausdorff and pathwise connected topological spaces. 
If $f:X \rightarrow  Y$ is a  proper local homeomorphism, then there exists $k\in \mathbb N$ such that for any $y\in Y$, $\card(f^{-1}(y))=k$ ($f$ is a ``$k$-fold covering''). Moreover, if $Y$ is simply connected, then $k=1$ i.e. $f$ is a homeomorphism.
\end{theorem}

Given the polyhedral domain  $\Omega\subset  \mathbb{R}^n$, we denote by    $\mathcal T$ be a triangulation of $\oo$. A subset $F$ of $\partial\Omega$ is called a \emph{facet}, or a $(n-1)$-dimensional face, of $\Omega$ if 
(i) $F$ is a polyhedral domain of an hyperplane $H$ (called the \emph{supporting hyperplane} of $F$) and 
(ii)  if $\sigma$ is a $(n-1)$-simplex in $\partial\Omega$ such that $\sigma\cap F$ contains a $(n-2)$-simplex, then either $\sigma\subset F$ or $\sigma \not\subset H$.
For any $0\leq k\leq n-2$, $G\subset \partial \Omega$ is called a 
$k$-dimensional face of $\Omega$ if it is a facet of a $(k+1)$-dimensional face.
Let $F_1,\dots,F_\ell$ be facets of a polyhedral domain $\Omega$ and $H_1,\dots,H_\ell$ be their supporting hyperplanes. We   say that $F_1,\dots,F_\ell$ are \emph{in general position} if $\dim (\bigcap_{i=1}^\ell H_i)=n-\ell$: {that is, the normals $\mathbf{n}_1,\ldots,\mathbf{n}_\ell$ to the hyperplanes are linearly-independent. Exploiting these definitions, we can show the following result: the proof is straightforward and is omitted.}

\begin{lemma}
\label{re:marque}
Let $\Omega$ be a bounded polyhedral domain of $\mathbb{R}^n$ and let $\mathcal T$ be a triangulation of $\oo$, and let  $\mathcal{F}_{k}$  be the  family of $k$-dimensional faces of $\mathcal{T}$
(cf. Definition \ref{def:polyhedron}). The following hold.
\begin{itemize}
\item
The boundary of $\oo$ is the  union of $(n-1)$-simplices of $\mathcal{F}_{n-1}$.
\item
 Since $\partial \Omega$ is a Lipschitz hyper-surface of $\mathbb{R}^n$ without boundary,  any $(n-2)$-simplex $\sigma\in \mathcal{F}_{n-2}$ contained in $\partial \Omega$ is the intersection locus of two $(n-1)$-simplices in  $\mathcal{F}_{n-1}$ that are contained in $\partial\Omega$.
\end{itemize}
\end{lemma}

{For three-dimensional domains, it is easy to verify that  vertices ($0$-faces) are the intersections of  $\ell\geq 3$ facets whose normals $\mathbf{n}_1,\ldots,\mathbf{n}_\ell$ span $\mathbb{R}^3$; similarly, 
the edges ($1$-faces) are the intersections of  $\ell'\geq 2$ facets whose normals $\mathbf{n}_1,\ldots,\mathbf{n}_{\ell'}$ define a two-dimensional subspace of $\mathbb{R}^3$.
 Lemma \ref{le:mme2} provides a generalization of this result to $n$-dimensional domains.
}

\begin{lemma}
\label{le:mme2}
Let $\Omega\subset  \mathbb{R}^n$ be a polyhedral domain. Any point of a $k$-dimensional face $G$ of $\Omega$ is contained in $n-k$ facets of $\Omega$ in general position. Furthermore, $G$ is contained in the intersection of the supporting hyperplanes of the facets.
\end{lemma}
\begin{proof}
We proceed by induction.
The property is obvious if $k=n-1$. 
Suppose that the result holds true for 
$k$+1-faces with 
$k+1\leq n-1$; we show below that the result holds for $k$-faces.

Let $\mathcal T$ be a triangulation of $\oo$ and let $G$ be a $k$-dimensional face of $\Omega$. 
Any point of $G$ lies in a $k$-simplex $\sigma\subset G$ that is contained in the boundary of a $(k+2)$-dimensional face $F$. More precisely, $\sigma$ lies in the intersection of two $(k+1)$-dimensional simplices $\tau$ and $\tau'$ of $\partial F$, cf. Lemma \ref{re:marque}.  Since $\sigma$ is in a $k$-dimensional face, it is not contained in the interior of a $k+1$-face and thus $\tau\cup\tau'$ is not contained in a facet of $F$. Hence $\sigma$ lies in the intersection of two non-parallel facets $E$ and $E'$ of $F$ and the supporting plane of $\sigma$ is the intersection locus of the supporting planes of $\tau$ and $\tau'$.

Since $E$ and $E'$ are $k+1$-dimensional faces of $\Omega$, any point in one of them belongs to the intersection of $n-k-1$ supporting hyperplanes in general position. Since $E$ and $E'$ are not contained in a $k+1$-dimensional affine subspace, for any point $\xi\in \sigma\subset E\cap E'$, it is possible to extract $n-k$ facets in general position among these $2(n-k-1)$ facets containing $\xi$.  The intersection of the corresponding supporting hyperplanes contains $\sigma$;  therefore, it  contains $G$.
\end{proof}

\subsection{Bijectivity of compositional maps}
We first present an extension of \cite[Proposition 2.1]{taddei2020registration}.
\begin{proposition}
\label{th:extended}
Let $k\geq0$ and $\Omega_0, \Omega_1, \dots \Omega_k \subset \mathbb{R}^n$ be bounded domains with continuous and connected boundary such that for any   $i,j \in \{1,\ldots,k\}$ 
(i)  $\oo_i\cap \oo_j=\emptyset$, if $i\neq j$; and 
(ii) $\overline{\Omega}_1, \ldots, 
\overline{\Omega}_k 
\subset \Omega_0$.

Let $\Omega=\Omega_0\setminus\bigcup_{i=1}^k\oo_i$ and $\Phi:\overline \Omega\mapsto \mathbb{R}^d$ be a locally injective continuous map such that $\forall i\in\{0,\dots,k\}$, $\Phi(\partial \Omega_i)\subseteq \partial\Omega_i$. 
If there exists $y\in \partial\Omega_0$ such that $\card(\Phi^{-1}(y)\cap \partial\Omega_0)=1$ or if $\partial \Omega_0$ is simply connected then  $\Phi(\oo)=\oo$ and $\Phi$ is a homeomorphism on its image. 
\end{proposition}
\begin{proof}
We split the proof in nine steps.

\noindent
{\bf 1. Existence of a tubular neighborhoods.}\\
According to Theorem \ref{th:brown},   for any $0\leq i\leq k$ there exists a 
tubular neighborhood
 $h_i$ from $\partial \Omega_i\times (-1,1)$ to a neighborhood of $\partial \Omega_i$ in $\mathbb{R}^n$ such that for any $\xi\in \partial \Omega_i$, $h_i(\xi,0)=\xi$. 
 \\[2mm]
\noindent
{\bf 2. The complements $C_i:=\mathbb{R}^n \setminus \oo_i$ are  pathwise connected and $\Omega$ is connected.}
\\
Indeed, consider two points $x$ and $y$ in $C_i$ and a path $\gamma:[0,1]\rightarrow \mathbb R^n$ connecting them. Let $\omega_i$ be a tubular neighborhood of $\partial \Omega_i$ and 
{$\omega_i^e=\omega_i\setminus \oo_i$}. From Theorem \ref{th:custom_theorem}, we deduce that if $\gamma$ does not intersect $\omega_i^e$ then it is contained in $C_i$. 
Otherwise,  there exist $0<t_0<t_1<1$ such that $\gamma([0,t_0])\subset C_i$, $\gamma([t_1,1]\subset C_i$ and $\gamma\{t_0,t_1\}\subset \omega_i^e$. 
Since, $\omega_i^e$ is homeomorphic to $\partial\Omega_i\times (0,1)$ and $\partial\Omega_i$ is connected, $\omega_i^e$ is connected. Therefore, it is possible to replace $\gamma|_{[t_0,t_1]}$ by a path in $\omega_i^e$ with the same endpoints. The path obtained is contained in $C_i$ and connects $x$ and $y$.
\\
Using the same reasoning, we can show that 
$\Omega=\Omega_0 \bigcap_{1\leq i\leq k} C_i$ is pathwise connected.
\\[2mm]
\noindent
{\bf 3. 
$\Phi\big|_{{\Omega}}$  (resp.  
$\Phi\big|_{\partial \Omega}$) is an  open map in $\mathbb{R}^n$ (resp. in $\partial \Omega$).
}\\ 
{We recall that the function $f:A\to B$ is an open map if $f(\omega)$ is open in $B$ for any open set $\omega$ in $A$.}
It suffices to show that, for any
  $y\in \Phi(\omega)$, there exists an open set $A\subset \Phi(\omega)$ that contains $y$. We denote by $\xi\in \omega$ an element of the pre-image of $y$.
  Since  $\Phi$ is locally injective, there exists $\eta>0$ such that $B_\eta(\xi)$, the open ball of radius $\eta$ centered at $\xi$, is contained in $\omega$ and  the restriction of 
  $\Phi$ to    $B_\eta(\xi)$ is injective. It follows from Theorem \ref{th:brouwer} that $A:=\Phi(B_\eta(\xi))$ is open.
 \\
{Since $\partial\Omega$ is locally the graph of a continuous function, it is locally homeomorphic to $\mathbb{R}^{n-1}$; therefore, we can apply the same argument to show that  $\Phi\big|_{\partial \Omega}$ is an open map.}
\\[2mm]
\noindent
 \textbf{4. $\mathbf{\partial  \Phi( \oo) \subseteq 
\Phi(\partial \Omega)  =\partial \Omega} $.}
\\
Let ${y \in \partial\, \Phi(\oo )}$. 
Since $\Phi$ is continuous and $\oo$ is closed, we find that $\Phi(\oo)$ is closed;
therefore, $\Phi^{-1}(y) $ is not empty.
Let $\xi\in \oo$ such that $\Phi(\xi)=y$. By contradiction, assume that $\xi\in \Omega$. Therefore, there exists $\eta>0$ such that the open  ball $B_\eta(\xi)\subset \Omega$. According to step \textbf{3.}, $\Phi( B_\eta(\xi)  )$ is an open set. Since,  $\Phi(B_\eta(\xi) ) \subset \Phi(\oo)$ and $y\in \Phi(B_\eta(\xi))$, $y\notin \partial \Phi(\oo)$. Contradiction. Hence $x\in \partial \Omega$ and $y\in \Phi(\partial \Omega)$.\\
{According to step {\textbf{3.}}, 
$\Phi(\partial \Omega_i)$ is open in $\partial\Omega_i$. Since $\partial \Omega_i$ is compact, it is also closed; as the boundary of $\Omega_i$ is connected,  Theorem \ref{th:connected_sets} implies that $\Phi(\partial \Omega_i)=\partial\Omega_i$ for $i=0,\ldots,k$.
}
\\[2mm]
\noindent
\textbf{5. $\mathbf{\Phi(\oo) \subseteq\oo_0}$ {\bf and} $\mathbf{\Phi^{-1}(\partial \Omega_0)=\partial \Omega_0}$.}
\\
By contradiction, let $z = \Phi(\xi) \in C_0=\mathbb{R}^d \setminus \oo_0$ for some $\xi \in \oo$.  
Since $\oo$ is compact and $\Phi$ is continuous, 
$\Phi(\oo)$ is bounded; since $C_0$ is unbounded (because $\Omega$ is bounded), there exists $y\in C_0 \setminus \Phi(\oo)$. Since $C_0$ is connected (cf. step \textbf{2.}), there exists a path 
 $\gamma \subset C$ that connects $z$ to $y$. 
 Recalling Theorem \ref{th:custom_theorem}, $\gamma$ should intersect $\partial \Phi(\oo)) \subseteq   \partial \Omega$, which contradicts the fact that $C_0\cap \partial \Omega = \emptyset$. Therefore, 
 $\Phi(\oo) \subseteq\oo_0$.
 \\
{Since   $\Phi(\Omega)$ is open (cf. step \textbf{3.}), we must have
$\Phi(\Omega)\cap \partial\Omega_0=\emptyset$ and thus
$\Phi^{-1}(\partial \Omega_0) :=\{x\in \overline{\Omega} : \Phi(x) \in \partial \Omega_0\} = \{x\in \partial \Omega : \Phi(x) \in \partial \Omega_0\}$.
Since   $\Phi(\partial \Omega_i) = \partial \Omega_i$ for $i=0,\ldots,k$  (cf. step \textbf{4.}), we find
$\Phi^{-1}(\partial \Omega_0) =\{x\in \partial \Omega_0 : \Phi(x) \in \partial \Omega_0\} =  \partial \Omega_0$.
}
\\[2mm]
\noindent
{\bf 6. There exists $\mathbf{y\in \partial\Omega_0}$ such that $\mathbf {\card\Phi^{-1}(y)=1}$.}
\\
From the latter step, we know that for any  $y\in \partial\Omega_0$ then $\Phi^{-1}(y)\subset \partial \Omega_0$ and thus
 $\card(\Phi^{-1}(y)) = \card(\Phi^{-1}(y)\cap \partial\Omega_0)$. 
 {Therefore, if there exists $y\in \partial\Omega_0$ such that  
 $\card(\Phi^{-1}(y)\cap \partial\Omega_0)=1$, we obtain the desired result.}
\\
Assume that $\partial\Omega_0$ is  simply-connected.
Since $\Phi\big|_{\partial\Omega_0}$ is a locally injective continuous map from $\partial \Omega_0$ onto itself, it follows from Theorem \ref{th:brouwer}, that it is a local homeomorphism. Since $\partial\Omega_0$ is compact this map is also proper. Moreover,  according to Theorem \ref{th:hadamard2}, $\Phi\big|_{\partial\Omega_0}$ is a homeomorphism. In particular, for any $y\in \partial\Omega_0$, $\card(\Phi^{-1}(y))=1$.
\smallskip

\noindent
 \textbf{7. Extension of $\mathbf \Phi$.}
\\ Let $h_i$ be a tubular neighborhood of $\partial \Omega_i$, as defined in step \textbf{1.}.
 We denote by $\omega_i$, the image of $h_i$ and by $\omega'_i$ the open set $h_i(\partial\Omega_i\times (-\frac 12,\frac 12))$. We choose $h_i$ such that $h_i(\partial \Omega_i\times [0,1))\subset \oo$ and  $\Phi(\omega'_i\cap \oo)\subset \omega_i$. \\
{We observe that for any $\xi\in \partial \Omega_i$ there exists $\eta>0$ such that $\Phi \big|_{B_\eta(\xi) \cap \overline{\Omega}}$ is injective and (cf. step \textbf{3.}) 
 $\Phi  ( B_\eta(\xi) \cap \partial \Omega_i )$ is an open neighborhood of $\Phi(\xi)$. We hence find 
\begin{equation}
\label{eq:tricky_identity_proof_difficult}
\forall \, \xi \in \partial \Omega_i, \; \exists \, \eta>0 \, : \, 
\Phi\left(  B_\eta(\xi) \cap \Omega \right) \cap \partial \Omega_i = \emptyset.
\end{equation}
 If \eqref{eq:tricky_identity_proof_difficult} is false, there exists a sequence $\{ \xi_n \}_n \subset \Omega$ that converges to $\xi$ such that $\Phi(\xi_n)\in \partial \Omega_i$ for all $n>0$. Since $\lim_{n\to \infty}\Phi(\xi_n) = \Phi(\xi)$ and {
 $\Phi({B_\eta(\xi) \cap \partial \Omega_i})$ is an open subset of $\partial \Omega_i$,}
  $\Phi(\xi_n) $ belongs to $\Phi\left(  B_\eta(\xi) \cap \partial \Omega_i  \right) $ for sufficiently large values of $n$: this contradicts the injectivity of $\Phi$ in $ B_\eta(\xi) \cap  \overline{\Omega}$. 
 }
 \\ 
Let  $\sigma_i:\omega_i\rightarrow \omega_i$ be the homeomorphism defined by $\sigma_i(h_i(\xi,t))=h_i(\xi,-t)$, for any $(\xi,t)\in \partial\Omega_i\times (-1,1)$. We define then $\Psi_i:\Omega\cup \omega'_i \rightarrow \mathbb{R}^n$ by $\Psi_i(x)=\Phi(x)$, if $x\in \oo$ and $\Psi_i(x)=\sigma_i(\Phi(\sigma_i(x))$ if $x \in \omega'_i\setminus \Omega$. {
The map $\Psi_i$ is continuous because it is the composition of continuous maps;
furthermore, it is locally injective  around any $\xi\notin \partial \Omega_i$,  since  $\Phi$ and $\sigma_i$ are both locally injective. \\
Let $\xi\in \partial \Omega_i$ and define $r>0$ such that $\Phi\big|_{B_r(\xi) \cap \overline{\Omega} }$  and 
$\Phi\big|_{\sigma_i(B_r(\xi) \setminus \overline{\Omega}  )}$
are injective. 
Let  
  $x\neq y\in B_r(\xi)$ belong to the same side of $\partial \Omega_i$. By construction, if $x,y\in \overline{\Omega}$, $\Phi(x)\neq \Phi(y)$ and thus $\Psi_i(x)\neq \Psi_i(y)$; similarly,
 if $x,y\notin \overline{\Omega}$, 
we have  $\sigma_i(x) \neq \sigma_i(y)$ and 
 $\sigma_i(x), \sigma_i(y) \in B_\epsilon(\xi)  \cap  \Omega$ for some $\epsilon>0$, and thus $\Psi_i(x) = \Phi(\sigma_i(x))\neq
 \Phi(\sigma_i(y)) =  \Psi_i(y)$. 
Let  $x \in \overline{\Omega}$ and $y \notin \overline{\Omega}$ and let $z=\sigma_i(y)\in \Omega$.
By contradiction, $\Psi_i(x)=\Psi_i(y)$, that is
 $\Phi(x)=\sigma_i(\Phi(z))$. Therefore,  $\Phi(x)$ and  $\Phi(z)$ are on opposite sides of $\partial \Omega_i$: since $\Phi$ is continuous, there exists $p \in {B_r(\xi) \cap} \Omega $ such that $\Phi(p) \in \partial \Omega_i$, which contradicts \eqref{eq:tricky_identity_proof_difficult}.
 In conclusion, we have shown that for any distinct 
$x,y\in B_r(\xi)$, 
$\Psi_i(x) \neq \Psi_i(y)$: the function $\Psi_i$ is hence locally injective.
}\\
 We can apply the same argument
 for all connected components of $\partial \Omega$:  we hence  obtain a continuous and locally injective  extension of $\Phi$ to an open neighborhood of $\oo$, that we call $\Psi$. Theorem \ref{th:brouwer} implies that  $\Psi$ is a local homeomorphism. 
\smallskip

\noindent
{\bf 8. The restriction of $\mathbf\Phi$ to $\mathbf{\Phi^{-1}(\oo)}$ is a homeomorphism onto $\mathbf\oo$.}
\\
Let $\omega$ be a connected component of $\Phi^{-1}(\overline \Omega)$. We consider $\Phi\big|_{\overline \omega}:\overline \omega\rightarrow \oo$. There exists an open set $\omega^+\supset\overline \omega$ on which $\Psi$ is well-defined and such that $\Psi(\omega^+\setminus\omega)\cap \Omega=\emptyset$. If $U$ is an open set of $\omega^+$, then, according to step \textbf{7.} $\Psi(U)$ is open. Moreover $\Psi(U)\cap \oo=\Psi(U\cap\overline\omega)=\Phi(U\cap\overline\omega)$. Hence, $\Phi$ sends open sets of $\overline \omega$ to open sets of $\oo$. Therefore, $\Phi\big|_{\overline \omega}:\overline \omega\rightarrow \oo$ is a local homeomorphism. 
\\
It follows that $\Phi(\overline\omega)$ is open and closed (it is compact) in $\oo$. According to step \textbf{2.}, $\oo$ is connected; 
since $\Phi(\overline\omega)$ is not empty, we can apply
 Theorem \ref{th:connected_sets} to conclude that $\Phi(\overline\omega)=\oo$. In particular, $\Phi^{-1}(\partial \Omega_0)\subset \overline \omega$. If $\omega'$ is another  connected component of 
 $\Phi^{-1}(\overline{\Omega})$ then $\partial\Omega_0\subset \overline \omega\cap \overline \omega'$. 
 Consequently, $\omega\cap \omega'\neq \emptyset$  and thus $\omega=\omega'$ i.e. $\Phi^{-1}(\overline{\Omega})$ is connected.
\\
We can now apply Theorem \ref{th:hadamard2} to $\Phi\big|_{\Phi^{-1}(\oo)}:\Phi^{-1}(\oo)\rightarrow \oo$. Since there exists $y\in \partial\Omega_0$ such that $\card(\Phi^{-1}(y))=1$, it is a homeomorphism. 
\smallskip

\noindent
 {\bf 9. Conclusion}
\\
We should   prove that
$\Phi^{-1}(\overline{\Omega})  =
 \{ x\in \oo : \Phi(x) \in \oo\} = 
  \overline{\Omega}$.
\\
Since $\Phi(\partial \Omega_i) = \partial \Omega_i$ for all $i$, we must have that 
$\partial \Omega \subset 
\Phi^{-1}(\overline{\Omega})$.
Since $\Phi\big|_{\Phi^{-1}(\overline{\Omega})}$ is a homeomorphism,
$\Phi(\partial  (
\Phi^{-1}(\overline{\Omega})
)
   ) = \partial \Phi( 
\Phi^{-1}(\overline{\Omega})   
    )
= \partial \Omega$.
Recalling step \textbf{ 4.}, 
we also have $\Phi(\partial \Omega) = \partial \Omega$.
In conclusion, we proved that 
$\Phi\big|_{\Phi^{-1}(\overline{\Omega})}$
 is injective and
 $\partial \Omega \subset 
\Phi^{-1}(\overline{\Omega})$; therefore, we have 
$\partial \Omega = \partial  (
\Phi^{-1}(\overline{\Omega})
)$.\\
Any $\xi\in \overline{\Omega}$ and $\xi'\in \Phi^{-1}(\overline{\Omega})$ can be connected by a path $\gamma$ in $\Omega$. The path  $\gamma$ does not cut $\partial\Omega=\partial \Phi^{-1}(\oo)$; therefore, according to Theorem \ref{th:custom_theorem}, $\xi\in \Phi^{-1}(\oo)$. 
\end{proof}

We comment on the assumptions of Proposition \ref{th:extended}.
The assumption $\card (\Phi^{-1}(y) \cap \partial \Omega) = 1$ for some $y\in \partial \Omega_0$ or alternatively the assumption that $\partial \Omega_0$ is simply connected is needed in step \textbf{6.}   to exclude examples where $\Phi(\Omega)=\Omega$ but $\Phi$ is not injective (but a $k$-fold covering with $k\neq 1$).
To provide a concrete example, consider the solid torus $\Omega=\Omega_0$ in $\mathbb{R}^3\simeq \mathbb C\times \mathbb{R}$, obtained by rotating a vertical disc around the axis $\{0\}\times\mathbb{R}$. The map $\Phi: (z,t)\mapsto (z^2,t)$ maps $\Omega$ in itself, it  is a local diffeomorphism (away from the axis $\{0\}\times\mathbb R$), but it is not injective:
note that  
$\card (\Phi^{-1}(y) \cap \partial \Omega) = 2$ for all $y\in \partial \Omega_0$ and 
the boundary of the torus is not simply connected.
Next, we  show a  characterization of
 connected bounded domains with continuous boundary.

\begin{lemma}
\label{le:mur}
Let $\Omega\subset \mathbb{R}^n$ be a connected bounded domain with continuous boundary. There exists  bounded domains with continuous and connected boundary  $\Omega_0, \Omega_1, ..$
$ \Omega_k \subset \mathbb{R}^n$, such that, for any  $i,j \in \{1,\ldots,k\}$,
(i)  $\oo_i\cap \oo_j=\emptyset$, if $i\neq j$; 
(ii) $\overline{\Omega}_i\subset \Omega_0$; 
(iii)  $\Omega=\Omega_0\setminus \bigcup_{i=1}^k\oo_i$.
\end{lemma}
\begin{proof}
Since $\Omega$ is bounded and its boundary is continuous, $\partial \Omega$ has a finite number of connected components, that we denote $\Sigma_0, \dots \Sigma_k$. The complement of each of them has two connected components, cf. step \textbf{2.} of Proposition \ref{th:extended}. One of them contains $\Omega$. We denote by $\Omega_i$ the bounded component of $\mathbb{R}^n\setminus \Sigma_i$. Since $\Omega$ is bounded, we can assume, without loss of generality, that $\Omega\subset \Omega_0$: since $\Omega$ is connected, we hence find that $\Sigma_i\subset \Omega_0$ for $i=1,\ldots,k$.
This  implies that  $\oo_i \subset \Omega_0$  and that  $\oo_1,\ldots, \oo_k$ are pairwise disjoint. Therefore, we obtain  $\Omega=\Omega_0\setminus \bigcup_{i=1}^k\oo_i$. 
\end{proof}

We need one additional Lemma.
\begin{lemma}
\label{th:tedious_lemma2}
Let $\Omega$ be a polyhedral domain  and 
\emph{$\Phi= \texttt{id} + \varphi$} with 
$\varphi \in \mathcal{U}_0$ satisfy $\min_{x \in \overline{\Omega}} J(\Phi ) > 0$. If $F$ is a facet of $\Omega$, then $\Phi(F)=F$. Furthermore,  $\Phi\big|_{\partial \Omega}$ is a homeomorphism onto $\partial \Omega$.  
\end{lemma}
\begin{proof}
Let   $F$ be a facet of $\Omega$ and let  $H$ be its supporting hyperplane with  normal $\mathbf{n}$;  since, for any $\xi\in F$, $\varphi(\xi)\cdot \mathbf{n}=0$,  we have $\Phi(F)\subset H$. Furthermore,   since,  $\min_{x \in \overline{\Omega}} J(\Phi ) > 0$, $\Phi$ is locally injective.

We proceed by induction on the dimension of the faces. Let $v$ be a vertex, i.e. a $0$-dimensional face.  
According to Lemma \ref{le:mme2}, there exists $n$ facets $F_1,\dots F_{n}$ with normals $\mathbf{n}_1,\dots \mathbf{n}_{n}$  such that $v=\bigcap_i F_i$ and $\mathrm{span}(\mathbf{n}_1,\dots \mathbf{n}_{n})=\mathbb{R}^n$.  Since, $\varphi(v)$ is orthogonal to each $\mathbf{n}_i$, $\varphi(v)=0$ and $\Phi(v)=v$.

Let $e$ be an edge of $\Omega$. It is an interval bounded by two vertices $v_-$ and $v_+$. Since any locally-injective continuous map in an interval is injective,   $\Phi\big|_e$ is injective. Exploiting  the intermediate value theorem, we must have $\Phi(e)=e$. 

Suppose that 
$\Phi(F)=F$ and $\Phi\big|_F$ is injective for any
 $k$-dimensional face of $\Omega$ with 
$1\leq k\leq n-2$.
Let $G$ be a $(k+1)$-dimensional face of $\Omega$ and let $P$ be its $(k+1)$-dimensional supporting plane.   Using again Lemma \ref{le:mme2},  any $\xi\in G$ belongs to  $\bigcap_{i=1}^{n-k-1} F'_i$, where the $F'_i$ are  facets, with supporting hyperplanes $H_i$. It follows from the above observation that $\Phi(\xi)\in \bigcap_{i=1}^{n-k-1} \Phi(F'_i)\subset \bigcap_{i=1}^{n-k-1} H_i=P$. Therefore $\Phi(G)\subset P$.
Since  $G$ is a polyhedral domain (cf. Lemma \ref{re:marque}), it is a connected bounded domain of $P$. Moreover, $\Phi\big|_G$ is continuous and locally injective and $\Phi\big|_{\partial G}$ is injective (it is a union of $k$-dimensional faces).   Therefore,
applying the result of  Proposition \ref{th:extended}  to $G\subset P$ and $\Phi\big|_G$, we have  $\Phi(G)=G$ and $\Phi\big|_G$ is injective. 

Hence, for any facet $F$ of $\Omega$, we obtain that $\Phi(F)=F$ and  $\Phi\big|_F$ is injective. 
Since the boundary of $\Omega$ is the union of facets, we find the desired result.
\end{proof}

We can now show the final result.
\begin{proof}
(\emph{Proposition \ref{th:bijectivity}}).
Since $\Omega$ is a  polyhedral domain, according to  Lemma \ref{le:mur},    there exists $k+1$ bounded
domains with continuous and connected boundary $\Omega_0,\dots,\Omega_k$ such that $\oo_i\cap \oo_j=\emptyset$ if $0<i<j$ and $\Omega=\Omega_0\setminus \bigcup_{i=1}^k\oo_i$. Moreover,  
since $\Omega$ has Lipschitz boundary, if $\Phi= \mathrm{id} + \varphi \in C^1(\oo,\mathbb R^n)$ and $|\mathrm {det}(\nabla \Phi(\xi)|>0$ for all $\xi\in \overline U$, we can resort to the inverse function theorem to show that 
$\Phi$ is  locally injective.
Finally,  Lemma \ref{th:tedious_lemma2} implies that $\Phi\big|_{\partial \Omega_0}$ is injective. 
Therefore, exploiting  Proposition \ref{th:extended},
we conclude that   $\Phi$ is a $C^1$ bijection and thus it is a diffeomorphism.
\end{proof} 

\subsection{Approximation properties of compositional maps}

\begin{proof}
(Proposition  \ref{th:approximation_DB}) According to Proposition 
\ref{th:approx_flows}, there exists a time-de- pendent vector field  $v_\Phi \in \mathcal{V}_0(\Omega)$  such that $\Phi=F[v_\Phi]:=X(v,1)$, where  
$(\xi,t) \mapsto  X(\xi, t)$ is the corresponding flow --- note that
$v_\Phi$ is tangent to the facets of $\Omega$. 
Recalling Lemmas \ref{re:marque} and \ref{le:mme2},
we have that any point of a $k$-dimensional face of $\Omega$ is in the intersection locus of $n-k$ facets in general position;  therefore,  $v_\Phi$ is tangent to any face of any dimension.  

Next, we consider a facet $F$ of $\Omega$. It is a polyhedral domain of an affine hyperplane $H$ with normal denoted $\mathbf n$:  the restriction of $v_\Phi$ to $F$ is hence a time-dependent vector field tangent to $\partial F$ (i.e. an element of $\mathcal{V}_0(F)$).  
Since the solution
to \eqref{eq:flow_diffeomorphisms} is unique for any Lipschitz-continuous velocity,
the flow of $v_\Phi\big|_F$ is the restriction of $X$ to $F$, that is  $\Phi\big|_F=F[v_\Phi\big|_F]$. 
Since the  facet $F$ is a Lipschitz domain of $H$ and $v_\Phi$ is tangent to its boundary,
we can apply
 Proposition \ref{th:bijectivity_flows} to conclude  that $\Phi(F)=F\subset H$. Since the function $\xi\mapsto \xi\cdot \mathbf n$ is  constant on $H$, for any $\xi \in F$, we find
$$
\varphi(\xi)\cdot \mathbf{n} 
=
\left(
\Phi(\xi)-\xi
\right)\cdot \mathbf{n} 
=\Phi(\xi) \cdot \mathbf{n}-\xi\cdot \mathbf{n}=0,
$$
which is the desired result.
\end{proof}

\appendix

\section{Existence of minimizers and performance of truncated bases}
\label{sec:easy_results}

\begin{proof}
(Proposition \ref{th:existence_minimizers})
Since $E$ is bounded from below, there exists $m^\star \in \mathbb{R}$ such that
$\inf_{\mathbf{a}\in \mathbb{R}^M} E(\mathbf{a})  \geq m^\star$. Given $R=E(0)  +  \lambda \mathfrak{f}_{\rm pen}(0) - m^\star$, we find that
$$
E(\mathbf{a}) + R \geq E(0) + 
\lambda \mathfrak{f}_{\rm pen}(0) = 
 \mathfrak{f}^{\rm obj}(0),
 \quad
 \forall \, \mathbf{a} \in \mathbb{R}^M.
$$

Since $\lim_{\| \mathbf{a} \|_2 \to \infty}
\mathfrak{f}_{\rm pen}(\mathbf{a}) = +\infty$, there exists $C>0$ such that
$\mathfrak{f}_{\rm pen}(\mathbf{a}) > R$ for any $\mathbf{a}$ that satisfies 
$\| \mathbf{a} \|_2 \geq C$. 
Given 
$\mathcal{B}_C(0):= \left\{ 
\mathbf{a} \in  \mathbb{R}^M \;
 {\rm s.t.} 
  \; \| \mathbf{a} \|_2 <  C \right\}$, 
  we hence find
\begin{equation}
\label{eq:lame_result}
 \mathfrak{f}^{\rm obj}(\mathbf{a})
 \geq 
E(\mathbf{a}) + R \geq E(0) + 
\lambda \mathfrak{f}_{\rm pen}(0) = 
 \mathfrak{f}^{\rm obj}(0),
 \;\;
 \forall \, \mathbf{a} \notin 
 \overline{
\mathcal{B}_C(0)}.
\end{equation}
 Since 
 $\overline{\mathcal{B}_C(0)}$ is a compact set in $\mathbb{R}^M$, we can apply Weierstrass theorem to conclude that the function  admits a minimum 
$\mathbf{a}^\star$
 in  $\overline{\mathcal{B}_C(0)}$:  
recalling \eqref{eq:lame_result}, 
$\mathbf{a}^\star$ is also a minimum of 
$ \mathfrak{f}^{\rm obj}$
 in $\mathbb{R}^M$. 
\end{proof}

\begin{proof}
(Lemma \ref{lemma:eigen_justification})
We first notice that
\begin{equation}
\label{eq:immediate_observation}
f( \mathbf{u}^\star )
+
\xi \|  
\mathbf{u}^\star 
\|_{\mathbf{A}}^2
\leq
f
(  P_W \mathbf{u}^\star )
+
\xi \|  
 P_W \mathbf{u}^\star 
\|_{\mathbf{A}}^2.     
\end{equation}
The eigenvectors
$\{ \boldsymbol{\phi}_i  \}_{i=1}^N$
satisfy $ \boldsymbol{\phi}_i^\top 
\mathbf{A} 
\boldsymbol{\phi}_j  = \delta_{i,j} \lambda_i$ and 
$\boldsymbol{\phi}_i^\top 
\mathbf{M} 
\boldsymbol{\phi}_j = \delta_{i,j} $ for $i,j=1,\ldots,N$ and  $ \delta_{i,j}$ is the Kronecker delta; we thus find
$
 \|  
\mathbf{u}^\star 
\|_{\mathbf{A}}^2 =  
\| P_W  \mathbf{u}^\star  \|_{\mathbf{A}}^2 
+
\|  P_W^\perp  \mathbf{u}^\star \|_{\mathbf{A}}^2 
$. Therefore, 
substituting the latter in \eqref{eq:immediate_observation},
recalling that $f$ is Lipschitz continuous and that
$\| P_W  \mathbf{u}^\star  \|_{\mathbf{M}} \leq
\|    \mathbf{u}^\star  \|_{\mathbf{M}}$,
we find
$$
\xi \|  
 P_W^\perp \mathbf{u}^\star 
\|_{\mathbf{A}}^2
\leq
f
(  P_W \mathbf{u}^\star )
-
f( \mathbf{u}^\star )
\Rightarrow 
 \| P_W^\perp \mathbf{u}^\star \|_{\mathbf{A}}^2
\leq
\frac{L^\star}{\xi}
 \| P_W^\perp \mathbf{u}^\star \|_{\mathbf{M}}.
$$
We further notice that
$$
 \| P_W^\perp \mathbf{u}^\star \|_{\mathbf{M}}^2
 =
 \frac{1}{\lambda_{m+1}}
  \underbrace{
 \sum_{i=m+1}^N
\boldsymbol{\phi}_i^\top \mathbf{M} \mathbf{u}^\star
 \,
 \lambda_i}_{=  \| P_W^\perp \mathbf{u}^\star \|_{\mathbf{A}}^2   }  \underbrace{\lambda_i^{-1}  \lambda_{m+1}}_{\leq 1}
 \leq
 \frac{1}{\lambda_{m+1}}
 \| P_W^\perp \mathbf{u}^\star \|_{\mathbf{A}}^2;
$$
we hence conclude that 
$ \| P_W^\perp \mathbf{u}^\star \|_{\mathbf{A}} 
\leq
\frac{L^\star}{\xi \sqrt{\lambda}_{m+1}}$ and 
$ \| P_W^\perp \mathbf{u}^\star \|_{\mathbf{M}} 
\leq
\frac{L^\star}{\xi  \lambda_{m+1}}$, which is \eqref{eq:eigen_justification_1}.

Finally, we have
$ \mathfrak{E} \leq \mathfrak{E}_m$ and also
$$
\mathfrak{E}_m - \mathfrak{E}
=
\underbrace{
f(  P_W \mathbf{u}^\star )
-
f
(  \mathbf{u}^\star )}_{
\leq L^\star 
\| P_W^\perp  \mathbf{u}^\star \|_{\mathbf{M}}    }
+
\underbrace{ 
\xi  \big( 
\|  
 P_W \mathbf{u}^\star 
\|_{\mathbf{A}}^2 
-
\|  
\mathbf{u}^\star 
\|_{\mathbf{A}}^2  \big)  }_{ \leq 0}
\leq
L^\star 
\| P_W^\perp  \mathbf{u}^\star \|_{\mathbf{M}}
\overbrace{\leq}^{\eqref{eq:eigen_justification_1}}
\dfrac{(L^\star)^2}{\xi \lambda_{m+1}},
$$
which is \eqref{eq:eigen_justification_2}.
\end{proof}

\section{Properties of vector flows}
\label{sec:vector_flows}

\subsection{Mathematical background}
Next theorem addresses the problem of extensions in $C^k$ spaces.
The result is due to Whitney
\cite{whitney1934functions}; further analyses and generalizations are provided in \cite[Chapter 2.5]{brudnyi2011methods} (cf. \cite[Theorem 2.64]{brudnyi2011methods}).
The result exploits the  definition of quasiconvex domains, which is provided below: we notice that any Lipschitz bounded domain is quasiconvex.
Extension theorems have been the subject of extensive research in Analysis; we refer to the monographs 
\cite{adams2003sobolev,brudnyi2011methods} for a thorough review of the subject. In particular, we refer to
\cite[Theorem 5.24]{adams2003sobolev} for an alternative result that addresses  the extension problem in  Sobolev spaces.  

\begin{Definition}
\label{def:quasiconvex_domain}
A bounded domain \( \Omega \subset \mathbb{R}^n \) is called \emph{quasi-convex} if there exists a constant \( C > 0 \) such that for every pair of points \( x, y \in \Omega \), there exists a rectifiable curve \( \gamma \subset \Omega \) joining \( x \) to \( y \) with the property that
\[
\text{length}(\gamma) \leq C \|x - y\|_2.
\]
The smallest such constant \( C \) is called the \emph{quasi-convexity constant} of the domain.
\end{Definition}

\begin{theorem}
\label{th:extension_whitney}
Let $\Omega$ be a quasiconvex domain and let $\varphi \in C^1(\Omega)$. 
If $\nabla \varphi$ can be defined on the boundary $\partial \Omega$ so that it is continuous in $\overline{\Omega}$, then  there exists a $C^1$ extension of $\varphi$ in $\mathbb{R}^n$. Furthermore, the extension operator is linear and bounded.
\end{theorem}

\subsection{Proofs}

\begin{proof} (Proposition \ref{th:bijectivity_flows})
Since $\Omega$ is quasi-convex, we have that $\Omega\times (0,1)$ is also quasi-convex; we can hence apply 
 the strong extension theorem \ref{th:extension_whitney} to prove the existence of a $C^1$ extension 
$V$  of 
the velocity $v$ to  $\Omega \times (0,1)$. Since the extension operator is continuous, $t\mapsto V(\cdot, t)$ is also continuous. Without loss of generality, we assume that $V$ is compactly supported in $\mathbb{R}^n \times [0,1]$; that is, it vanishes for $\|x\|_2 \geq R$ for some constant $R>0$ and $t\in [0,1]$.

We define the flow $\Psi$ of $V$, 
$\Psi=F[V]$.
Since $V$ is uniformly continuous  and its support is contained in the bounded ball $B_R(0)$, there exist $C_0,C_1>0$ such that $\|V(\xi,t) \|_2 \leq C_0 + C_1 \| \xi \|_2$ for all $\xi\in \mathbb{R}^n$ and $t\in [0,1]$. Therefore, 
recalling the Cauchy-Lipschitz theorem,
the solution to the ODE in \eqref{eq:flow_diffeomorphisms} exists and is unique, for all initial conditions $\xi\in \mathbb{R}^n$.  By differentiating \eqref{eq:flow_diffeomorphisms}, we obtain a formal  expression for $\nabla \Psi$
$$
\nabla \Psi(\xi) =\nabla  X(\xi,t=1),
\;\;
{\rm where} \;\;
\left\{
\begin{array}{ll}
\frac{\partial \nabla X}{\partial t}(\xi,t) = \nabla_x V( X(\xi,t) , t)  \nabla X(\xi,t) & t\in (0,1] \\[3mm]
\nabla  X (\xi,0) = \mathbbm{1} & \\
\end{array}
\right. 
$$
Since $\nabla_x V$ is uniformly continuous, the latter ODE is also well-posed.
Finally, 
since the field $V$ satisfies the hypotheses of 
 \cite[Theorem 8.7]{younes2010shapes}, we find that $\Psi$ is  a diffeomorphism in $\mathbb{R}^n$.

Provided  that $X(\overline{\Omega},t) =  \overline{\Omega}$ for all   $t\in [0,1]$,  the restriction of $\Psi$ to $\Omega$ is equal to $\Phi$ --- that is, 
$\Phi= \Psi\big|_{\Omega} $ ---  and hence that $\Phi$ is a bijection in $\overline{\Omega}$.

We now prove that 
$X(\overline{\Omega},t) =  \overline{\Omega}$ for all   $t\in [0,1]$ (see also
\cite[Theorem 4.4]{santambrogio2015optimal}).
Towards this end, we introduce
  $\rho_0 \in \mathcal{P}_2(\mathbb{R}^n) \cap C^\infty(\mathbb{R}^n)$ and we define the density of the  pushforward measure
$\rho(x,t) = \rho_0(\xi) \left( {\rm det} \big( \nabla X(\xi,t)  \big) \right)^{-1}$ with $x=X(\xi,t)$. 
Given $\phi\in C_{\rm c}^\infty(\mathbb{R}^n)$ 
---
$ C_{\rm c}^\infty(\mathbb{R}^n)$ is the set of 
$C^\infty(\mathbb{R}^n)$ functions with compact support --- 
and introducing the notation $ J(\xi, t): = {\rm det} \big( \nabla X(\xi,t)  \big)$,
 we find that\footnote{Note that here we need to consider integration over $\mathbb{R}^n$ because we do not know if $X(\cdot, t)$ is a diffeomorphism in $\Omega$.}
$$
\begin{array}{l}
\displaystyle{
\int_{\mathbb{R}^n}
\partial_t \rho(x,t) \phi(x) \, dx
=
\frac{d}{dt}
\int_{\mathbb{R}^n}
 \rho(x,t) \phi(x) \, dx
   }
\\[3mm]
\displaystyle{
\overset{\rm (CV)}{=}
\frac{d}{dt}
\int_{\mathbb{R}^n}
 \phi(X(\xi,t)) 
 \rho_0(\xi) \left( J(\xi, t)  \big) \right)^{-1} 
  J(\xi, t)     \, d\xi
  =
\frac{d}{dt}
\int_{\mathbb{R}^n}
 \phi(X(\xi,t)) 
 \rho_0(\xi)  \, d\xi
 }
\\[3mm]
\displaystyle{
 =
\int_{\mathbb{R}^n}
\frac{\partial}{\partial t}
 \phi(X(\xi,t)) 
 \rho_0(\xi)  \, d\xi  
  =
\int_{\mathbb{R}^n}
 \nabla_x \phi(X(\xi,t)) \cdot  \partial_t X(\xi,t)
 \rho_0(\xi)  \, d\xi
}
\\[3mm]
\displaystyle{
   =
\int_{\mathbb{R}^n}
 \nabla_x \phi(X(\xi,t)) \cdot  v( X(\xi,t), t)
 \rho_0(\xi)  \, d\xi
}
\\[3mm]
\displaystyle{  
  =
\int_{\mathbb{R}^n}
 \nabla_x \phi(X(\xi,t)) \cdot  v( X(\xi,t), t)
 \rho_0(\xi, t)   
 \left( J(\xi, t)   \right)^{-1} 
J(\xi, t) \, d\xi
 }
 \\[3mm]
 \displaystyle{
=
\int_{\mathbb{R}^n}
 \nabla_x \phi(X(\xi,t)) \cdot  v( X(\xi,t), t)
 \rho(X(\xi, t)   , t)
J(\xi, t) \, d\xi
 }
 \\[3mm]
 \displaystyle{
\overset{\rm (CV)}{=}
\int_{\mathbb{R}^n}
 \nabla_x \phi(x) \cdot  v( x, t)
 \rho(x  , t)
 \, dx
 }
 \\
\end{array}
$$
 where 
$\overset{\rm (CV)}{=}$ refers to the application of the change of variable formula. Provided that $v$ and $X$  (and thus $\rho$) are sufficiently regular, we find
$$
\int_{\mathbb{R}^n}
\left( 
\partial_t \rho(\cdot ,t)  
\, + \nabla_x \cdot \left( v( \cdot, t)  \rho(\cdot  , t) \right)
\, \right) 
   \phi   \, 
dx
= 0
\qquad
\forall \, \phi \in C_{\rm c}^\infty(\mathbb{R}^n),
\quad
\forall \, t\in [0,1].
$$
Since $C_{\rm c}^\infty(\mathbb{R}^n)$ is dense in $L^2(\Omega)$, we also find
$$
\int_{\Omega}
\left( 
\partial_t \rho(\cdot ,t)  
\, + \nabla_x \cdot \left( v( \cdot, t)  \rho(\cdot  , t) \right)
\, \right) 
   \phi   \, 
dx
= 0
\qquad
\forall \, \phi \in L^2(\Omega),
\quad
\forall \, t\in [0,1].
$$
Then, if we plug $\phi \equiv 1$ and we apply the Gauss theorem, we find
$$
\frac{d}{dt}
\int_\Omega
  \rho(x,t)     \, dx
  =
  \int_\Omega
  \partial_t \rho(x ,t)   \, dx
=
-
  \int_{\partial \Omega}
  \partial_t   \rho(x, t)  
\underbrace{ v(x ,t)   \cdot \mathbf{n}(x) }_{=0}  
 \, dx
  = 0
$$
and thus
\begin{equation}
\label{eq:santambrogio_explained}
\int_\Omega
  \rho(x,t)     \, dx
=  
\int_\Omega
  \rho_0(\xi)     \, d\xi,
  \quad
  \forall \, \rho_0\in \mathcal{P}_2(\mathbb{R}^n)\cap C^\infty(\mathbb{R}^n).
\end{equation}

The identity
\eqref{eq:santambrogio_explained}
 can be readily extended to the general case  through a density argument. We can indeed introduce  a sequence $\{v_j \}_{j=1}^\infty \subset C^\infty$ such that $\lim_j \|v - v_j  \|_{C^1(\Omega\times (0,1))} = 0$ and define $X_j$ as the flow of $v_j$, and $\rho_j$ as the density of the corresponding pushforward measure. Since $v_j \in C^\infty$, we can readily show that also $X_j \in C^\infty$  (cf. 
 \cite[Theorem 8.9]{younes2010shapes}) 
 and thus $\rho_j  \in C^\infty$.
Exploiting  \eqref{eq:continuity_data}, we can show that 
$X_j$ converges to $X$ 
in $L^\infty(\Omega)$;
furthermore,
using a similar argument\footnote{
If we define $W=\nabla X$ and $W_j = \nabla X_j$, we find that the difference $e_j=W-W_j$ satisfies (we omit dependence on $\xi$):
$$
\left\{
\begin{array}{ll}
\partial_t e_j = \nabla v  W - 
\nabla v_j  W_j =
\nabla v_j e  +(\nabla v - 
\nabla v_j)  W
& t\in (0,1] \\
e_j(0) = 0 & \\
\end{array}
\right.
$$
Recalling  Lemma \ref{th:tedious_lemma_odes},
 Eq. \eqref{eq:state_transition},
 and the fact that any converging sequence is bounded,
 we readily find that 
 $$
 \max_{\xi \in \overline{\Omega}}
 \| W(\xi,t) -   W_j(\xi,t) \|_{2}
 \leq
C \| \nabla v(\cdot, t)  - \nabla v_j(\cdot, t) \|_{L^\infty(\Omega\times (0,1))},
 $$
 which implies convergence of the sequence $\{W_j \}_j$ to $W$ in $L^\infty$.
}
we can also show that
 $\nabla X_j$ converges to $\nabla  X$ in $L^\infty(\Omega)$.
In conclusion, we find that  $\lim_j \|\rho(\cdot,t) -\rho_j(\cdot,t)  \|_{L^2(\Omega)} = 0$ for all $t\in [0,1]$, which implies that $\rho$ satisfies
\eqref{eq:santambrogio_explained}.

By contradiction, we assume that $X(\bar{\xi},t)\notin \Omega$ for some $\bar{\xi} \in \Omega$ and $t\in (0,1]$. Since 
$X(\cdot,t)$ is continuous and is a global diffeomorphism, there exists an open  neighborhood $A=\mathcal{B}_r(\bar{\xi}) \cap \Omega$ such that $X(\xi,t)\notin \Omega$ for all $\xi \in A$. Consider now $\rho_0:\Omega\to \mathbb{R}_+$ such that
(i) $\rho_0(\xi) \geq 0$, (ii) $\rho_0 \in C^\infty(\Omega)$,  
(iii)   $\rho_0(\xi) = 0 $ for all $\xi \notin A$, and 
(iv) $  \int_\Omega
  \rho_0(\xi)     \, d\xi=1$. 
Exploiting the definition of the 
pushforward density, 
since $\Omega \cap X(A, t)  = \emptyset$,
we find
$$
\int_\Omega
  \rho(x,t)     \, dx
  =
 \int_{\Omega \cap X(A, t) }
  \rho(x,t)     \, dx 
  = 0 \neq 1 = 
  \int_\Omega
  \rho_0(\xi)     \, d\xi,
$$ 
which contradicts \eqref{eq:santambrogio_explained}.
\end{proof}

\begin{proof} (Proposition \ref{th:approx_flows})
Since $\Phi \in {\rm Diff}_0(\overline{\Omega})$, 
by definition
there exists an isotopy $X: \overline{\Omega} \times [0,1] \to \overline{\Omega}$ such that
(i) $X(\xi, t=0)=\xi$ for all $\xi \in \overline{\Omega}$, (ii) $X(\xi,t=1)=\Phi(\xi)$, 
(iii)   $X(\cdot, t) \in {\rm Diff}_0(\overline{\Omega})$  for all $t\in [0,1]$.
Then, we can define the vector field $v$ such that
$v(X(\xi,t) , t):=\frac{\partial X}{\partial t}(\xi,t)$. By construction, $t \mapsto X(\xi, t)$ is the unique solution to  the ODE system \eqref{eq:flow_diffeomorphisms} for all $\xi\in \overline{\Omega}$. 

By contradiction, suppose that 
$v(x , \bar{t})\cdot \mathbf{n}(x) \neq 0$ for some $x\in \partial \Omega$ where the boundary $\partial \Omega$  is locally differentiable (and thus the normal $\mathbf{n}$ is properly defined) and $ \bar{t}\in (0,1]$.
Since $X(\cdot, \bar{t}) \in {\rm Diff}_0(\overline{\Omega})$, we must have 
$X(\partial \Omega, \bar{t}) = \partial \Omega$ and thus there exists a unique $\xi\in \partial \Omega$ such that $X(\xi, \bar{t}) = x$.
Next, we introduce the distance function 
$$
t \in [0,1] \mapsto d(\xi,t):={\rm dist}(X(\xi,t), \partial \Omega) ,
$$
associated with the trajectory 
$t \in [0,1] \mapsto X(\xi,t)$.
Since $X(\cdot, t) \in  {\rm Diff}_0(\overline{\Omega})$ for all $t\in [0,1]$, we must have
 that $d(\xi,t)\equiv 0$.
On the other hand, recalling that $\nabla_x {\rm dist}(x, \partial \Omega) = \mathbf{n}(x)$  for all $x\in \partial \Omega$ where $\partial \Omega$ is locally differentiable, we  find that
$$
 d(\xi, \bar{t}+ \Delta t)
=
\underbrace{ d(\xi, \bar{t})}_{=0}
+
\underbrace{\frac{\partial d}{ \partial t}(\xi, t)}_{= \partial_t X(\xi, \bar{t}) \cdot  \nabla_x {\rm dist}(x, \partial \Omega)  }
\Delta t
+ 
o(\Delta t)
=
v(x , \bar{t}) \cdot \mathbf{n}(x)
\Delta t
+ 
o(\Delta t),
$$
for all $|\Delta  \,  t| \ll 1$. This 
 implies that $t\mapsto  d(\xi, t)$
 does not vanish in a neighborhood of $t=\bar{t}$. Contradiction.
\end{proof}

\begin{proof}
(Proposition \ref{th:evolution_jacobian})
We observe that  \eqref{eq:flow_diffeomorphisms}  is a system of   ODEs that depends on the initial condition $\xi\in \Omega$ and on the parameter $\mathbf{a}\in \mathbb{R}^M$. 
By differentiating \eqref{eq:flow_diffeomorphisms} with respect to $\xi\in \Omega$, we find that $t \mapsto  \nabla X (\xi, t)$ satisfies
\begin{equation}
\label{eq:ODE_gradient}
\left\{
\begin{array}{ll}
\frac{\partial \nabla X}{\partial t}(\xi,t) = \nabla_x v( X(\xi,t) , t)  \nabla X(\xi,t) & t\in (0,1], \\[3mm]
\nabla  X (\xi,0) = \mathbbm{1}. & \\
\end{array}
\right. 
\end{equation} 

Recalling the Jacobi's formula for the derivative of the determinant, we find
$$
\begin{array}{rl}
\displaystyle{
\frac{\partial J}{\partial t} (\xi, t)
=}
&
\displaystyle{
 J (\xi, t)  {\rm trace} \left(  \left(
\nabla X(\xi, t )
\right)^{-1}
\frac{\partial \nabla X }{\partial t} (\xi, t )
\right)
}
\\[3mm]
=
&
\displaystyle{
=J (\xi, t)  {\rm trace} \left(  \left(
\nabla X(\xi, t )
\right)^{-1}
\nabla_x v \left( X(\xi, t ), t \right)
\,
\nabla X(\xi, t )
\right).
}
\\
\end{array}
$$
Since ${\rm trace}(A B) = {\rm trace}( B A) $ for any pair of compatible square matrices $A$ and $B$, we hence find
$$
\frac{\partial J}{\partial t} (\xi, t)
\, = \, 
J (\xi, t)  {\rm trace} 
\left(  
\nabla_x v \left( X(\xi, t ), t \right)
\right)
\, = \,
J (\xi, t)  \; 
\left(
\nabla_x \cdot  v \left( X(\xi, t ), t \right)
\right),
$$
which is the desired result.
\end{proof}

\begin{proof}
(\emph{Proposition \ref{th:approximation_flows}}).
We denote by $X$ the flow associated with $v$ and by $Y$ the flow associated with $w$.
Given $\xi\in \Omega$, we define $e(t):= \| X(\xi, t) - Y(\xi, t) \|_2$. Then, we have that
$$
\begin{array}{l}
\displaystyle{e(t) = 
\big\|
\int_0^t
\left(
\partial_t X(\xi, s)
-
\partial_t Y(\xi, s)
\right)
\, ds
\big\|_2
=
\big\|
\int_0^t
\left(
v(X(\xi, s),s)
-
w(Y(\xi, s),s)
\right)
\, ds
\big\|_2
}
\\[3mm]
\leq
\displaystyle{
\int_0^t
\big\|
v(X(\xi, s),s)
-
w(Y(\xi, s),s)
\big\|_2
\, ds
}
\\[3mm]
\leq
\displaystyle{
\int_0^t
\big\|
v(X(\xi, s),s)
-
v(Y(\xi, s),s)
\big\|_2
\, ds
+
\int_0^t
\big\|
v(Y(\xi, s),s)
-
w(Y(\xi, s),s)
\big\|_2
\, ds
}
\\[3mm]
\end{array}
$$
Exploiting the definition of Lipschitz constant and introducing $\delta(t):= 
\sup_{\xi \in \Omega}
\big\|
v(\xi,s)
$
$
-
w(\xi,s)
\big\|_2$, we obtain
$$
e(t) \leq
\int_0^t L e(s) \, ds
+
\int_0^t \delta(s) \, ds.
$$
Then, applying Gronwall's inequality and
observing that
$\max_{s\in [0,1]} \delta(s) \leq 
\| v -w  \|_{L^\infty(\Omega \times (0,1))}
$, we find
$$
\begin{array}{rl}
e(t) \leq
&
\displaystyle{
 \int_0^t 
\delta(s) e^{L(t-s)} \, ds
\leq
\| v -w  \|_{L^\infty(\Omega \times (0,1))}
\int_0^t 
 e^{L(t-s)} \, ds
}
\\[3mm]
&
\displaystyle{
 =
 \| v -w  \|_{L^\infty(\Omega \times (0,1))}
 \frac{e^{L t} -  1}{L},
}
\\
\end{array}
$$
which is the desired result.
\end{proof}

\subsubsection{Computation of the gradient of the target function}
In view of the proof of 
Propositions \ref{th:derivative_computations_VB} and
\ref{th:adjoint}, we recall the following Lemma.

\begin{lemma}
 \label{th:tedious_lemma_odes}
 (cf. \cite{baake2011peano})
 Consider the system of ODEs:
 \begin{equation}
 \label{eq:ODE_linear}
 \left\{
 \begin{array}{ll}
 \dot{x}(t) = A(t) x(t) + b(t) & t> 0, \\[3mm]
 x(0) = x_0, & \\
 \end{array}
 \right.
 \end{equation}
 where $A: \mathbb{R}_+ \to \mathbb{R}^{n\times n}$ and 
 $b: \mathbb{R}_+ \to \mathbb{R}^{n}$ are smooth functions.
 There exists a \emph{state transition matrix} $\Phi: \mathbb{R}_+ \times \mathbb{R}_+ \to  \mathbb{R}^{n\times n}$
which depends only on $A$  
  such that the solution to \eqref{eq:ODE_linear} is given by 
 \begin{equation}
 \label{eq:ODE_linear_sol}
x(t)  = \Phi(t,s) x(s) + \int_s^t \Phi(t,\tau) b(\tau) \, d\tau,
\quad
\forall \; t,s\in \mathbb{R}_+,  t>s.
 \end{equation}
 \end{lemma}

\begin{proof}
(Proposition    \ref{th:derivative_computations_VB}).
If we differentiate 
\eqref{eq:flow_diffeomorphisms} with respect to $a_i$, we find
\begin{equation}
\label{eq:dXda}
\left\{
\begin{array}{ll}
\displaystyle{
\frac{\partial }{\partial t}
\frac{\partial X}{\partial a_i}(\xi,t)
\, = \, 
\nabla_x  v \left( X(\xi,t),t \right) 
\frac{\partial X}{\partial a_i}(\xi,t)  \, + \,
\frac{\partial v}{\partial a_i} \left( X(\xi,t),t \right) 
}
& t\in (0,1) \\[3mm]
\displaystyle{
\frac{\partial X}{\partial a_i}(\xi,0) = 0},
&
\\
\end{array}
\right.
\end{equation}
for $i=1,\ldots,M$.
We notice that \eqref{eq:dXda} and \eqref{eq:ODE_gradient} read as parametric ODE systems of the form \eqref{eq:ODE_linear} with the same matrix $A(t) = \nabla_x  v \left( X(\xi,t),t \right) $; therefore, they share the same state transition matrix $\Phi$. 

Exploiting Lemma \ref{th:tedious_lemma_odes},  we hence find that 
\begin{equation}
\label{eq:state_transition}
\Phi( t, \tau; \xi) 
\, = \, 
\nabla X(\xi, t)  \left(  \nabla X(\xi, \tau)   \right)^{-1}.
\end{equation}
If we substitute \eqref{eq:state_transition} in 
  \eqref{eq:ODE_linear_sol},  we finally obtain
$$
\frac{\partial X}   {\partial a_i}(\xi, t)
=
\nabla X(\xi, t) 
\int_0^t \, 
\left(
\nabla X(\xi, \tau) 
\right)^{-1}
\phi_i \left( X(\xi, \tau), \tau \right) \, d\tau,
\quad
{\rm for} \;\; i=1,\ldots,M,
$$
which is \eqref{eq:derivative_VB}.
\end{proof}

\begin{proof}
(Proposition    \ref{th:adjoint}).
The proofs of \eqref{eq:adjoint_method_distributed} and
\eqref{eq:adjoint_method_pointwise} exploit the very same argument: we hence only rigorously prove 
\eqref{eq:adjoint_method_pointwise}.
Below, we omit dependence on $\mathbf{a}$ to shorten notation.
First,
given $k\in \{1,\ldots,M \}$,
 we rewrite \eqref{eq:sensitivity_method_b} as
\begin{equation}
\label{eq:adjoint_proof1}
\frac{\partial E}{\partial a_k} 
\, = \, 
\sum_{i=1}^{N_0} \,\Psi_i \cdot 
\frac{\partial \texttt{N} }{\partial a_k} (\xi_i) 
=
\sum_{i=1}^{N_0} \,\Psi_i \cdot 
\frac{\partial X_i}{\partial a_k} (1) 
\end{equation}
where $X_i(t)=X(\xi_i,t)$ and
$\Psi_i:=
 X_i(1) - 
 \left( \sum_{j=1}^{N_1} P_{i,j}  y_j  \right)$ for $i=1,\ldots,N_0$. Exploiting the definition of the adjoint, we find
$$
\int_0^1 \left(
\frac{d \Lambda_i}{d t}(t) \, + \,  \nabla_x v \left( X_i(t),t   \right) ^\top \, \Lambda_i(t)
\right)
\cdot \frac{\partial X_i}{\partial a_k}(t) 
\, dt
= 0,
\quad
{\rm for} \; i=1,\ldots,N_0.
$$
If we integrate by part the first term and we exploit the final  condition of 
\eqref{eq:adjoint_method_pointwise_b} and the fact that
$\frac{\partial X_i}{\partial a_k}(0) = 0 $, we find
\begin{equation}
\label{eq:adjoint_proof2}
\int_0^1 
\left(
\frac{d }{d t}
 \frac{\partial X_i}{\partial a_k}(t) 
 \, - \,  \nabla_x v \left( X_i(t),t   \right)   \, 
  \frac{\partial X_i}{\partial a_k}(t) 
\right)
\cdot
\Lambda_i(t)
\, dt
= 
\Psi_i \cdot  \frac{\partial X_i}{\partial a_k}(1),
\quad
{\rm for} \; i=1,\ldots,N_0.
\end{equation}
By comparing \eqref{eq:adjoint_proof1} with  \eqref{eq:adjoint_proof2}, we deduce
$$
\frac{\partial E}{\partial a_k} 
\, = \, 
\sum_{i=1}^{N_0}
\int_0^1 
\left(
\frac{d }{d t}
 \frac{\partial X_i}{\partial a_k}(t) 
 \, - \,  
  \nabla_x v 
 \left( X_i(t),t   \right)    \, 
  \frac{\partial X_i}{\partial a_k}(t) 
\right)
\cdot
\Lambda_i(t)
\, dt
$$ 
and then exploiting \eqref{eq:dXda} we obtain
$$
\frac{\partial E}{\partial a_k} 
\, = \, 
\int_0^1 
 \frac{\partial v}{\partial a_k} 
 \left( X_i(t),t   \right) 
\cdot
\Lambda_i(t)
\, dt,
$$ 
which is \eqref{eq:adjoint_method_pointwise_b}.
\end{proof}

\bibliographystyle{abbrv}
\bibliography{all_references}
  
\end{document}